\documentclass[11pt,fleqn]{article}
\usepackage{mathptmx}       % selects Times Roman as basic font
\usepackage{helvet}         % selects Helvetica as sans-serif font
\usepackage{courier}        % selects Courier as typewriter font
\usepackage{type1cm}        % activate if the above 3 fonts are
\usepackage{url}

\usepackage{makeidx}         % allows index generation
\usepackage{graphicx}        % standard LaTeX graphics tool
\usepackage{array,booktabs,calc}
\usepackage{multicol}        % used for the two-column index
 \usepackage{multirow}
\usepackage[bottom]{footmisc}% places footnotes at page bottom

\usepackage{cite}
\usepackage{fancyvrb}        % LaTeX tool for writing verbatim
\usepackage{amsmath}
\usepackage{amssymb}
\usepackage{algorithm,algorithmic}
\usepackage{url}

\usepackage{url}
\makeindex             % used for the subject index
\begin{document}
\title{Determining the essentially different partitions of all Japanese convex tangrams}
%\author{T.G.J.~Beelen\\[0.3cm]\today}
%\author{Theo G.J.~Beelen}
%\footnote{T.G.J.Beelen@tue.nl}
%}
%\institute{T.G.J.~Beelen \at Eindhoven University of Technology, the Netherlands,\\
%\email{{T.G.J.Beelen}@tue.nl}
%}
%
\author{T.G.J.~Beelen and T. Verhoeff (TU Eindhoven) \\[0.3cm]}
%\authorrunning{T.G.J.~Beelen\hspace{0.3cm}\today}
%
%\footnote{T.G.J.~Beelen \at Eindhoven University of Technology, Dep. Mathematics and Computer Science (CASA), P.O. Box 513, 5600 MB  Eindhoven, the Netherlands}
%\footnote{T.G.J.~Beelen,\email{T.G.J.Beelen@tue.nl}}
%
%
\maketitle
\tableofcontents
%-------------------------------------------------------------------------
\clearpage
\section{Abstract}
In this report we consider the set of the 16 possible convex tangrams that can be composed with the 7 so-called ``Sei Shonagon Chie no Ita'' (or Japanese) tans, see \cite{Tangram:ref_epstein_uehara_2014}. The set of these Japanese tans is slightly different from the well-known set of 7 Chinese tans with which 13 (out of those 16) convex tangrams can be formed.  
In \cite{Tangram:ref_beelen_2018}, \cite{Tangram:ref_pentoma_2017} the problem of determining all essentially different partitions of the 13 ``Chinese'' convex tangrams was investigated and solved. In this report we will address the same problem for the ``Japanese'' convex tangrams.
The approach to solve both problems is more or less analogous, but the ``Japanese'' problem is much harder than the ``Chinese'' one, since the number of ``Japanese'' solutions is much larger than the ``Chinese'' ones. In fact, only for a few ``Japanese''  tangram shapes their solutions can be found by a rigorous analysis supported by a large number of clarifying diagrams. The solutions for the remaining shapes have to be determined using a dedicated computer program.
Both approaches will be discussed here and all essentially different solutions with the ``Japanese'' tans are presented.
As far as we know all presented results are not yet published before.  \newline
\newline\newline
\mbox{\it{Keywords:}}\;\; tangram, partition, backtracking, visualization.   
\clearpage
\section{General Introduction}\label{sec: Gen_Intro}
\begin{figure}[thb]
\begin{center}
\includegraphics[scale = 0.4, keepaspectratio]{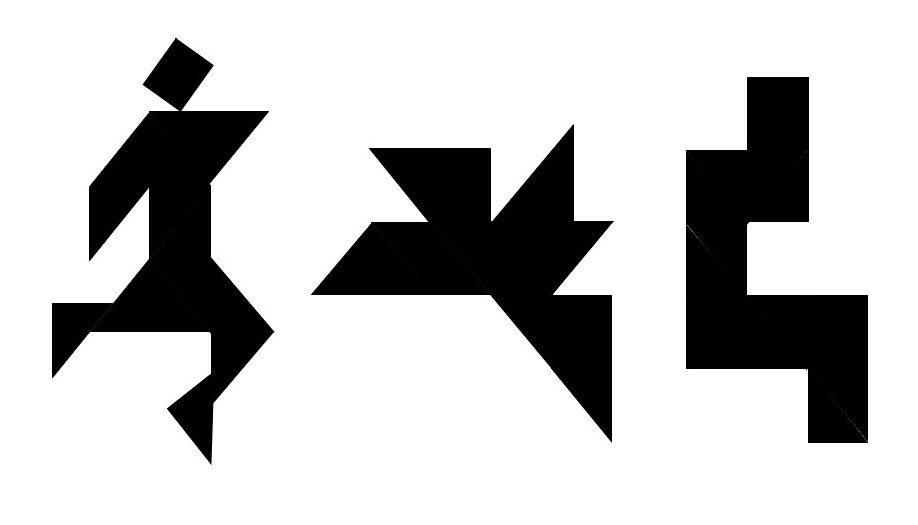}
\caption{A few typical tangram figures.}
\label{fig: Tangram_figures_org}		%	fig. 1
\end{center}
\end{figure}

The word ``tangram'' is reasonably well-known, especially in the context of trying to compose a given figurative picture using 7 geometrically shaped puzzle pieces called tans. A few typical examples of such pictures are given in Fig.~\ref{fig: Tangram_figures_org}. 
Furthermore, in the bilingual (German / Dutch) book ``Tangram, / Das alte chinesische Formenspiel / Het oude Chinese vormenspel '' by J. Elffers \cite{Tangram:ref_elffers_1976} 
over 1600 examples on tangram puzzles and their solutions can be found. See also website \cite{Tangram:ref_pentoma_2017}  \url{http://www.pentoma.de/}.
\begin{figure}[thb]
\begin{center}
\includegraphics[scale=0.4, keepaspectratio]{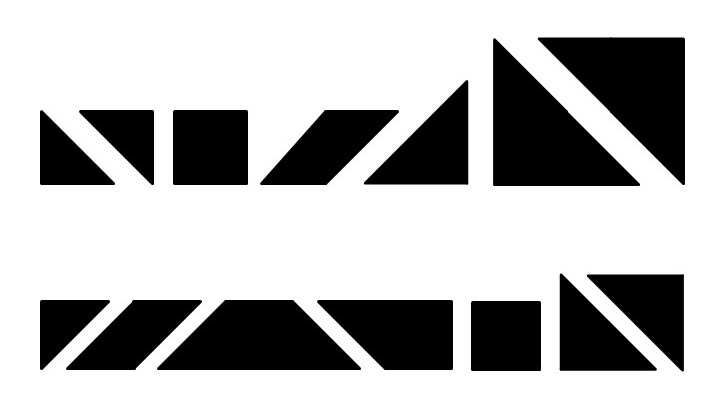}
\caption{The 7 individual Chinese and the 7 Japanese tangram pieces, in the top and bottom row, respectively.}
\label{fig: Chin_Jap_tans_1}		%fig. 2		
\end{center}
\end{figure}%

\begin{figure}[thb]
\begin{center}
\includegraphics[scale=0.6, keepaspectratio]{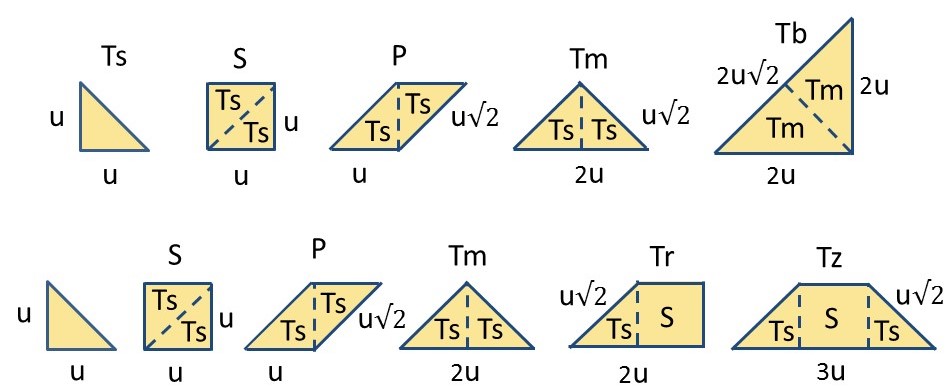}
\caption{The composition of the Chinese (top row) and Japanese tans (bottom).}
\label{fig: Chin_Jap_tans_2}		%fig. 3		
\end{center}
\end{figure}

\noindent
It should be noticed that next to the set of Chinese tans, there also exists a less known set of Japanese tans, in shape similar to the Chinese tans. Both the sets of Chinese and Japanese tans are shown in  Fig.~\ref{fig: Chin_Jap_tans_1}. 
Just as in case of the Chinese tans, each Japanese piece can be decomposed into one or more of the smallest triangular pieces. See Fig.~\ref{fig: Chin_Jap_tans_2}. 
\subsection{The relationships between the tans}
Let us consider Figs.~\ref{fig: Chin_Jap_tans_1} and ~\ref{fig: Chin_Jap_tans_2} in more detail. The set of Chinese tans consists of the following 7 pieces: one small triangle $Ts$, one square $S$, one pararallelogram $P$, two medium sized triangles $Tm$ and two big triangles $Tb$. \newline
The set of Japanese tans also consists of 7 pieces: one small triangle $Ts$, one square $S$, one pararallelogram $P$, two medium sized triangles $Tm$ (exactly as in the Chinese set), one rectagular trapezium $Tr$ and one isosceles trapezium $Tz$. Notice that the Chinese and Japanese tans $Ts,\;S,\;P$ and $Tm$ are identical.
\newline\newline
The relationship between the areas of the tans is indicated in Fig.~\ref{fig: Chin_Jap_tans_2}.
\begin{eqnarray}
	& & \hspace*{-1.5cm}\mbox{\rm{\bf{Relationships\; between\; the\; tans\; :}}}		\label{lab: Tan_relations} \\
	& & \hspace*{-0.8cm}\mbox{\rm{\;\;(i)\; S = 2\,Ts , \quad\; P = 2\,Ts ;}} \nonumber\\
		& & \hspace*{-0.8cm}\mbox{\rm{\;(ii)\; Tm = 2\,Ts, \quad 
		Tb = 2\,Tm = \,4\,Ts ; }} \nonumber\\
%extra for Jap. tans
		& & \hspace*{-0.8cm}\mbox{\rm{\;(iii)\; Tr = Ts + S = 3\,Ts, \quad 
		Tz = 2\,Ts + S = Tr + Ts = 4\, Ts ; }} \nonumber\\		
	& & \hspace*{-0.8cm}\mbox{\rm{(iv)\; Tangram with all 7 Chinese tans \;\, =\, 2\,Tb + Tm\, +\; 2\,Ts \;\, + \; S +\, P\, = \, 16\;Ts ;}} \nonumber\\
	& & \hspace*{-0.8cm}\mbox{\rm{\,(v)\; Tangram with all 7 Japanese tans \, =\, Tz + Tr + 2Tm\, + Ts + \, S +\, P\, = \, 16\;Ts ;}}	 \nonumber
\label{tan_relations}
\end{eqnarray}	
\subsection{Problem formulation}
For convenience, let us denote the above mentioned sets of Chinese and Japanese tans by $Tans\_C$ and $Tans\_J$, respectively. 
It is immediately clear that we can create a huge number of tangrams that can be composed using either the set $Tans\_C$ or $Tans\_J$. In 1942 two Chinese mathematicians Fu Traing Wang and Chuan-Chih Hsiung published a paper \cite{Tangram:ref_paper_1942} dealing with the problem of finding all $\mathit{convex}$ polygons in 2D. 
%
%\begin{figure}[thb]
%\begin{center}
%\includegraphics[scale=0.6, keepaspectratio]{All_convex_polygons}
%\caption{All 20 convex polygons. }
%\label{fig: All_convex_polygons}		
%\end{center}
%\end{figure}
%
\begin{figure}[thb]
\begin{center}
\includegraphics[scale=0.6, keepaspectratio]{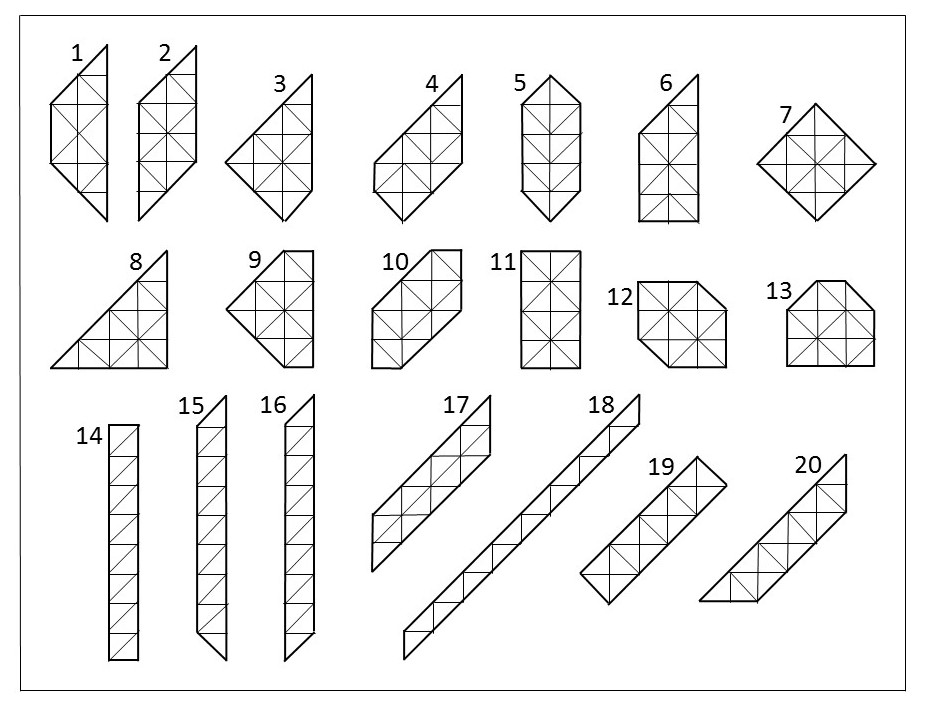}
\caption{All 20 convex polygons. }
\label{fig: grid_tangrams}		
\end{center}
\end{figure}
%\begin{figure}[t!hb]
%\begin{center}
%\hspace{-1.4cm}
%\includegraphics[scale=0.6, keepaspectratio]{All_Chin_Jap_polygons}
%\caption{Left: The 13 convex polygons that can be formed  by the Chinese tans. \newline Right: The 16 convex polygons that can be formed  by the Japanese tans. }
%\label{fig: All_Chin_Jap_polygons}		
%\end{center}
%\end{figure}
They showed that there exist precisely 20 convex polygons, see Fig.~\ref{fig: grid_tangrams}. Here the numbering is conform to that in \cite{Tangram:ref_paper_1942} as follows. The shapes marked by $*$ in their table correspond to the numbers 14-16, while the 13 unmarked  shapes as listed in their order are numbered 1-13 here. Clearly, these 13 shapes can be covered by the set $C\_tans$ and they will be denoted by the set $Poly13$.\newline
Notice that all 20 polygons can be formed by 16 rectangular isosceles triangles, as shown in Fig.~\ref{fig: grid_tangrams}.
\newline
In 2014 it was shown by Eli Fox-Epstein and Ryunhei Uehara \cite{Tangram:ref_epstein_uehara_2014} that even 3 more (so in total 16, but not more) convex polygons can be covered by using the set $J\_tans$, being the shapes 14-16 in Fig.~\ref{fig: grid_tangrams}. We will call the polygons 1-16 the set $Poly16$  (thus $Poly16$ $\supset$ $Poly13$). 
Moreover, it was also shown in \cite{Tangram:ref_epstein_uehara_2014} that more coverings are possible, but then different tan sets must be used.
\newline\newline
A mathematically interesting question on tangram covering is how many essentially different partitions for all polygons in $Poly13$ and $Poly16$ exist. In \cite{Tangram:ref_beelen_2017} this problem was extensively analyzed and solved for $Poly13$ with the set $C\_tans$.\newline
In this report we will address the covering problem for $Poly16$ with the set $J\_tans$.\newline\newline
$\mathbf{Terminology}:$ In the rest of this report we will use the words $\mathit{layout,\; partition,\; filling}$ and $\mathit{covering}$ as synonyms. 
\section{Finding all layouts of the convex polygons with the Japanese tans}
In this section we will give an overview of a few approaches we have chosen to solve the problem of finding all essentially different layouts with the Japanese tans (in short, to solve ``The problem'').
First we notice that solving this problem can be done by the well-known backtracking technique \cite{Tangram:ref_wikipedia_2017}, \cite{Tangram:ref_geeks_2017}, in a similar way as done in case of the Chinese tans (see \cite{Tangram:ref_beelen_2018}).\newline
However, after proceeding manually as in \cite{Tangram:ref_beelen_2018} for a few ``simple'' tangrams it became clear that this would be a very tedious job for almost all convex shapes. So, an alternative was to carry out a backtracking algorithm by a computer. We will address this approach in section \ref{sec: backtr_algo}.
In fact, we will discuss the following approaches to solve (partly or fully) ``The problem''.\newline\newline
(1)  A systematic analysis of all possible partitions for the full square tangram, i.e. tangram 7 in Fig.~\ref{fig: grid_tangrams}.
\newline\newline
(2) A systematic analysis of all possible partitions for the rectangular strips (see shapes 14 up to 16 in Fig.~\ref{fig: grid_tangrams}, which will be indicated by $J14,\; J15$ and $J16$ in the rest of this report). This analysis is done via backtracking and visualizing all possible trees.
%Notice that due to the length of the analyses of $J15$ and $J16$, they are placed in the Appendix.
Notice that due to the length of the analyses of $J15$ and $J16$, this is not included here, but can be found in \cite{Tangram:ref_beelen_verhoeff_2018}.
\newline\newline
(3) An alternative `combinatorial' approach for the strips $J14,\; J15$ and $J16$.
\newline\newline
(4) A global description of an algorithm to solve a so-called packing problem. We will explain how this algorithm can be used to solve ``The Problem''. Moreover, all solutions (partitions and their number) of all 16 convex shapes are included.
\clearpage
%
%					 The Square tangram
%
\section{The $Square$ tangram}
\begin{figure}[thb]
\begin{center}
\includegraphics[scale = 0.39]{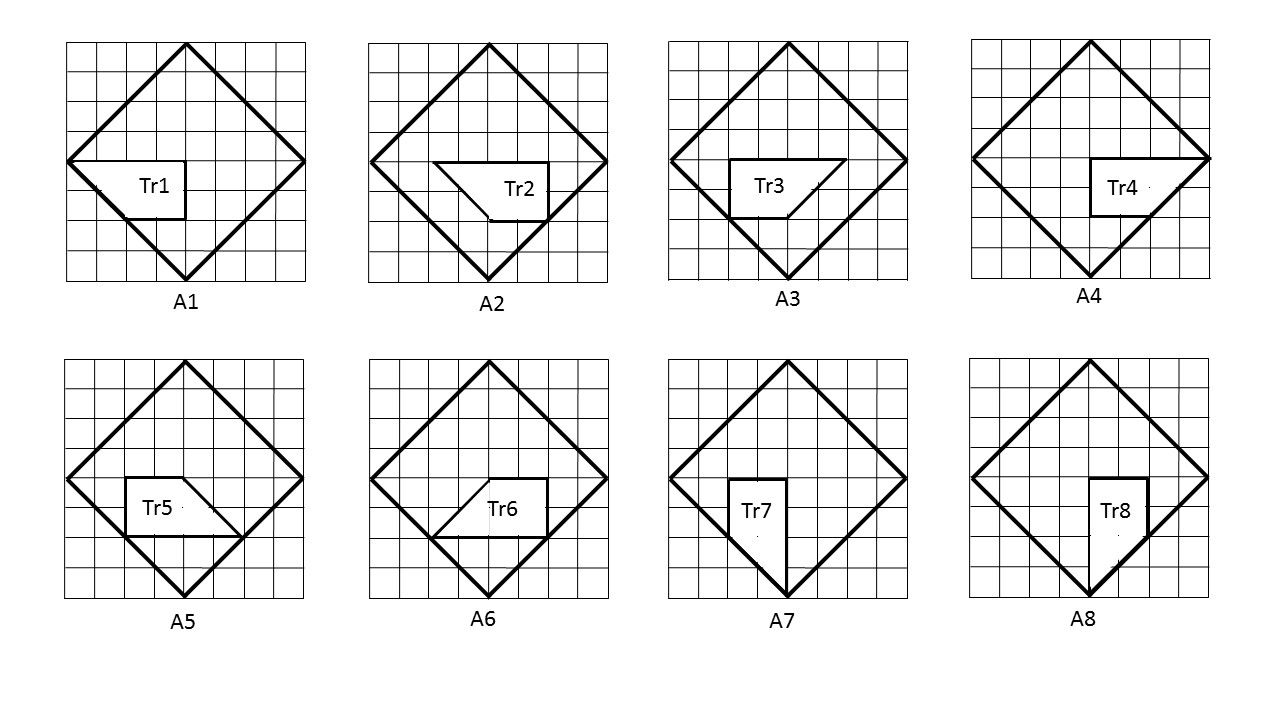}
\caption{All possible layouts for tan Tr in the lower half of  tangram $Square$}
\label{fig: Square_all_Tr}		%fig 6
\end{center}
\end{figure}

We start with adding the rectangular trapezium Tr to the empty tangram $Square$. Since $Square$ is symmetric w.r.t. both its horizontal and vertical centerline we can restrict ourselves to placing Tr in the lower half of $Square$. See Fig.~\ref{fig: Square_all_Tr}.
We will discuss several cases in the next sections.
\subsection{The case Tr1}
\begin{figure}[thb]
\begin{center}
\includegraphics[scale = 0.4]{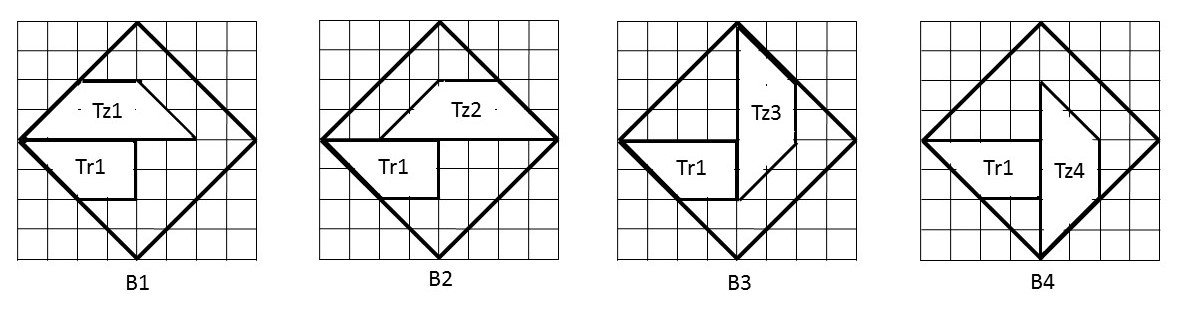}
\caption{All possible combinations when adding Tz to Tr1 in $Square$}
\label{fig: Square_Tr1_Tz}		%fig 7
\end{center}
\end{figure}

Let us first consider the case with Tr1, see Fig.~\ref{fig: Square_all_Tr}-A1. Clearly, we can add the isosceles trapezium Tz on 4 different positions as shown in Fig.~\ref{fig: Square_Tr1_Tz}.
\begin{figure}[thb]
\begin{center}
\includegraphics[scale = 0.4]{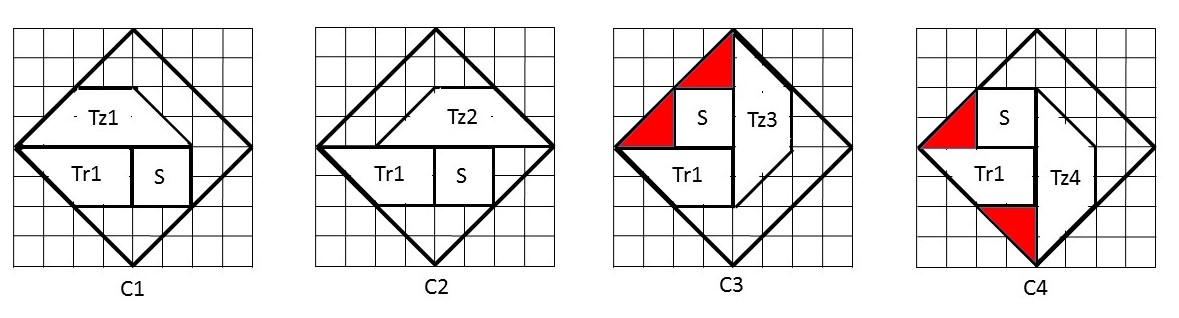}
\caption{All possible combinations when adding S to Tr1.Tz in $Square$}
\label{fig: Square_Tr1_Tz_S}		%fig 8
\end{center}
\end{figure}
Next we can add the square S to these 4 layouts. Clearly, we have only one possibility for S per layout. See Fig.~\ref{fig: Square_Tr1_Tz_S}. However, the layouts C3 and C4 are not feasible since 2 tans Ts are needed (but not available) for a full covering of the square. So we have
\begin{eqnarray}
	& & \hspace*{-0.7cm}\mbox{\rm\bf{Conclusion:}} 	
	\label{lab: Conclusion_feas_C1_C2} \\		%concl 
	& & \hspace*{-0.7cm}\mbox{\rm\bf{The $Square$ tangram has 2 potentially feasible layouts C1 and C2 with Tr1 and Tz}.} \nonumber
\end{eqnarray}
\subsection{The case Tr4}
\begin{figure}[thb]
\begin{center}
\includegraphics[scale = 0.4]{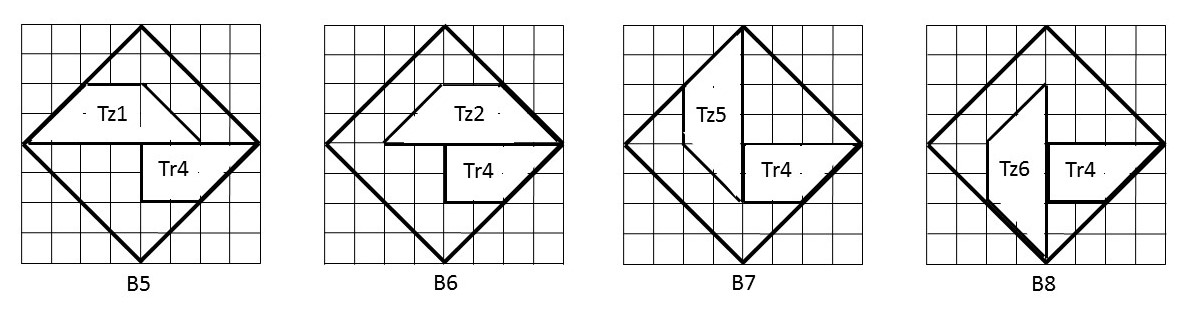}
\caption{All possible combinations when adding Tz to Tr4 in $Square$}
\label{fig: Square_Tr4_Tz}		%fig 9
\end{center}
\end{figure}
\begin{figure}[thb]
\begin{center}
\includegraphics[scale = 0.4]{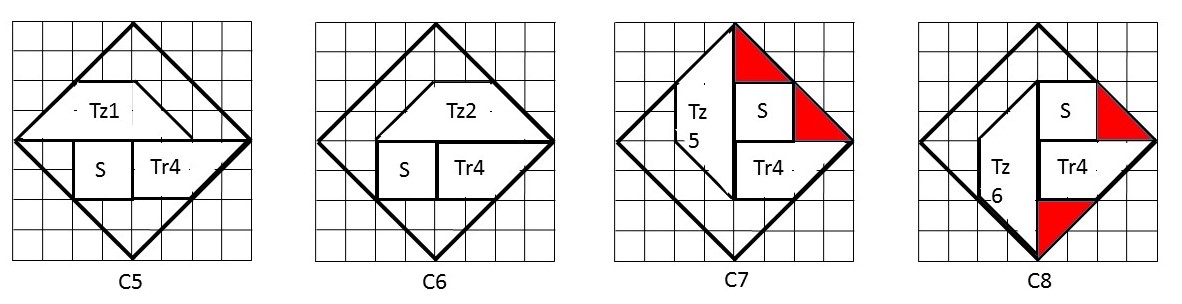}
\caption{All possible combinations when adding S to Tr4.Tz in $Square$}
\label{fig: Square_Tr4_Tz_S}		%fig 10
\end{center}
\end{figure}%

\noindent
In a similar way as in case Tr1 we can proceed with case Tr4. Indeed, here again we find 4 different combinations for Tr4 and Tz as shown in Fig.~\ref{fig: Square_Tr4_Tz}.
Next we can add S to these 4 layouts, see Fig.\ref{fig: Square_Tr4_Tz_S}. Similarly to case Tr1 the layouts C7 and C8 for Tr4 are not feasible. 
So we have 
\begin{eqnarray}
	& & \hspace*{-0.7cm}\mbox{\rm\bf{Conclusion:}} 	
	\label{lab: Conclusion_feas_C5_C6} \\		%concl 
	& & \hspace*{-0.7cm}\mbox{\rm\bf{The $Square$ tangram has 2 potentially feasible layouts C5 and C6 with Tr4 and Tz}.} \nonumber
\end{eqnarray}
When comparing the potentially feasible layouts in Figs.~\ref{fig: Square_Tr1_Tz} and \ref{fig: Square_Tr4_Tz} we see that 
the cases C1 and C6 as well as C2 and C5 are $\mathit{equivalent}$ (by symmetry w.r.t. the vertical centerline of $Square$). Thus, 
\begin{eqnarray}
	& & \hspace*{-0.7cm}\mbox{\rm\bf{Conclusion:}} 	
	\label{lab: Conclusion_study_C1_C2} \\		%concl 
	& & \hspace*{-0.7cm}\mbox{\rm\bf{We have to investigate the cases C1 and C2 further, but}} \nonumber\\
	& & \hspace*{-0.7cm}\mbox{\rm\bf{we can skip the cases C5 and C6 with Tr4}.} \nonumber
\end{eqnarray}
Before doing this, we first turn to the cases with Tr2, Tr3, Tr5 and Tr6. 
\subsection{The cases Tr2, Tr3, Tr5 and Tr6}
\begin{figure}[thb]
\begin{center}
\includegraphics[scale = 0.4]{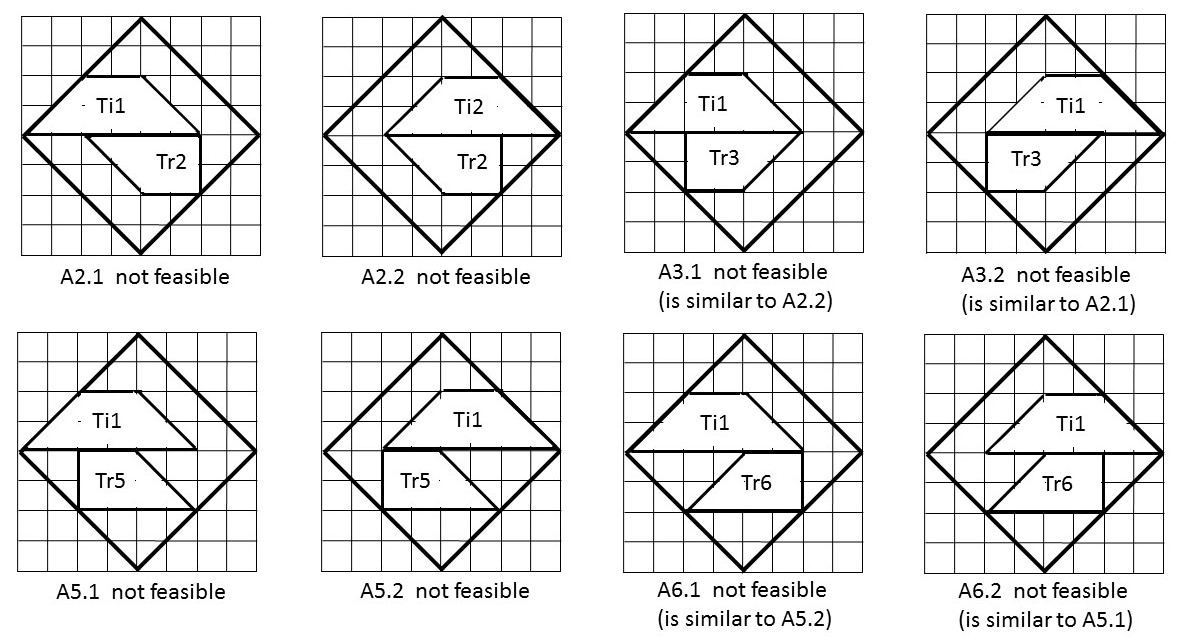}
\caption{All possible combinations when adding Tz to Tr2, Tr3, Tr5 and Tr6 in $Square$}
\label{fig: Square_Tr_nonfeas}		%fig 11
\end{center}
\end{figure}%

\noindent
In Fig.~\ref{fig: Square_Tr_nonfeas} all possible cases with Tr2, Tr3, Tr5 and Tr6 are shown. It is easily seen that in all cases we cannot add the square tan S anymore. So, have 
\begin{eqnarray}
	& & \hspace*{-0.7cm}\mbox{\rm\bf{Conclusion:}} 	
	\label{lab: Tr2356_not_feas} \\		%concl 
	& & \hspace*{-0.7cm}\mbox{\rm\bf{The cases Tr2, Tr3, Tr5 and Tr6 are not feasible.}} \nonumber	
\end{eqnarray}
\newline
Finally, we want to study the cases with Tr7 and Tr8.
\subsection{The cases Tr7 and Tr8}
\begin{figure}[thb]
\begin{center}
\includegraphics[scale = 0.4]{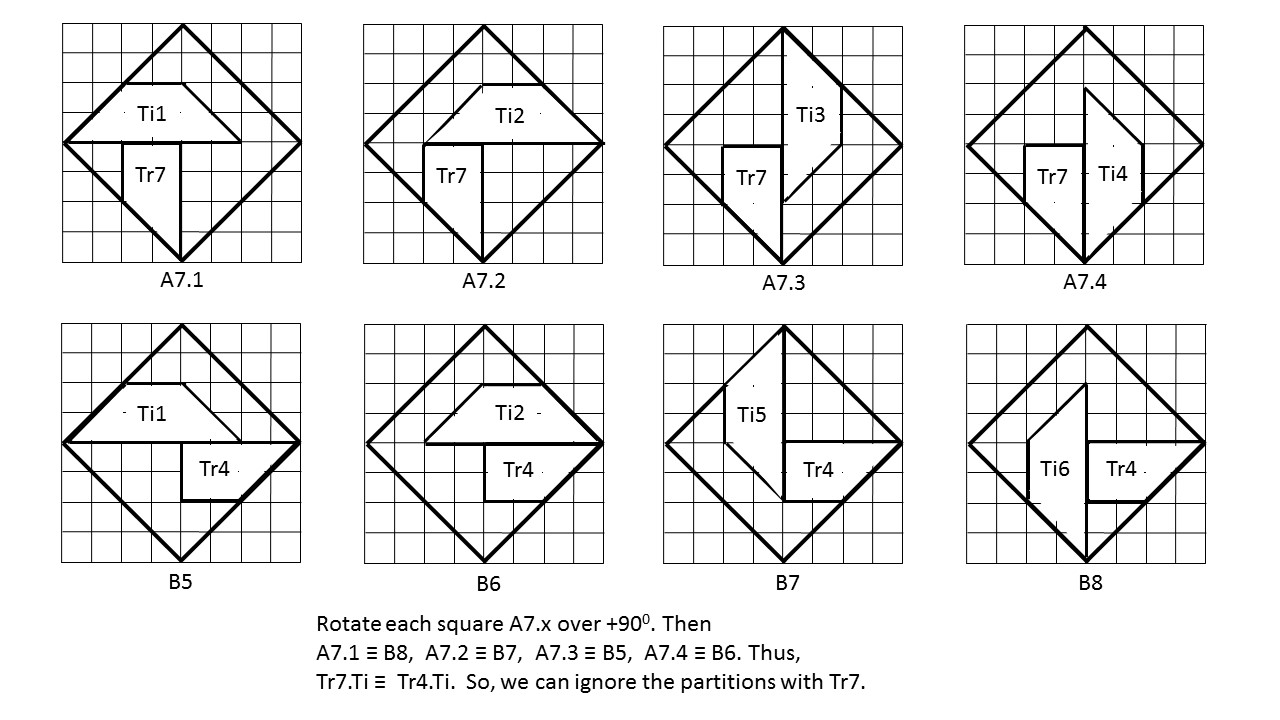}
\caption{All possible combinations when adding Tz to Tr7 as well as to Tr4 in $Square$}
\label{fig: Square_Tr7}		%fig 12
\end{center}
\end{figure}
\begin{figure}[thb]
\begin{center}
\includegraphics[scale = 0.4]{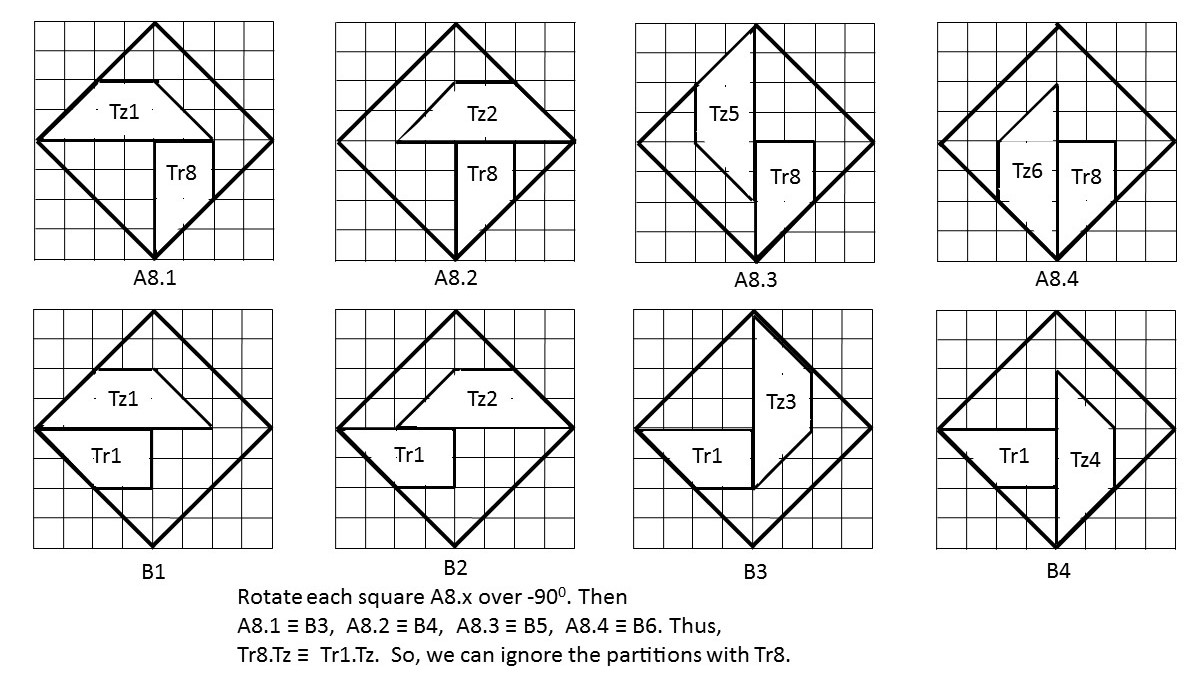}
\caption{All possible combinations when adding Tz to Tr8 as well as to Tr1 in $Square$}
\label{fig: Square_Tr8}		%fig 13
\end{center}
\end{figure}

In Fig.~\ref{fig: Square_Tr7} all 4 possible cases for Tr7 in $Square$ are shown. 
Furthermore, the 4 possible cases for Tr4 are repeated here (recall Fig.\ref{fig: Square_Tr4_Tz}). It can easily be seen that by rotating each configuration $Square\_Tr7.Tz$ over $90^0$ we find one of the configurations with $Square\_Tr4.Tz$. We call such a pair $\mathit{equivalent}$  and this will be denoted by the symbol $\equiv$.
Specifically, in Fig.~\ref{fig: Square_Tr7} we have $A7.1 \equiv B8,\; A7.2 \equiv B7,\; A7.3 \equiv B6$ and $A7.4 \equiv B5$. 
Hence, we have
\begin{eqnarray}
	& & \hspace*{-0.7cm}\mbox{\rm\bf{Conclusion:}} 	
	\label{lab: Conclusion_Ignore_Tr7} \\		%concl 
	& & \hspace*{-0.7cm}\mbox{\rm\bf{We can ignore all partitions with Tr7 for further study (and continue with Tr4).}} \nonumber
\end{eqnarray}
\newline
Let us now consider case $Square\_Tr8$.\newline
In Fig.~\ref{fig: Square_Tr8} all 4 possible cases for $Square\_Tr8$ are shown. Moreover, for easy comparison we also recall the 4 possible cases for Tr1 from Fig.\ref{fig: Square_Tr1_Tz}. It can easily be seen that by rotating each configuration $Square\_Tr8.Tz$ over $90^0$ we find one of the configurations with $Square\_Tr1.Tz$. 
Specifically, we have $A8.1 \equiv B3,\; A8.2 \equiv B4,\; A8.3 \equiv B2$ and $A8.4 \equiv B1$. 
Hence, we have
\begin{eqnarray}
	& & \hspace*{-0.7cm}\mbox{\rm\bf{Final Conclusion:}} 	
	\label{lab: Conclusion_Ignore_Tr8} \\		%concl 
	& & \hspace*{-0.7cm}\mbox{\rm\bf{We can ignore all partitions with Tr8 for further study (and continue with Tr1).}} \nonumber	
\end{eqnarray}
\newline
Combining all conclusions (\ref{lab: Conclusion_feas_C1_C2}) up to (\ref{lab: Conclusion_Ignore_Tr8}) above we see that we only need to investigate the layouts C1 and C2 in Fig.~\ref{fig: Square_Tr1_Tz_S}		for yes/no feasibility.
\clearpage
\begin{figure}[thb!]
\begin{center}
\includegraphics[scale = 0.4]{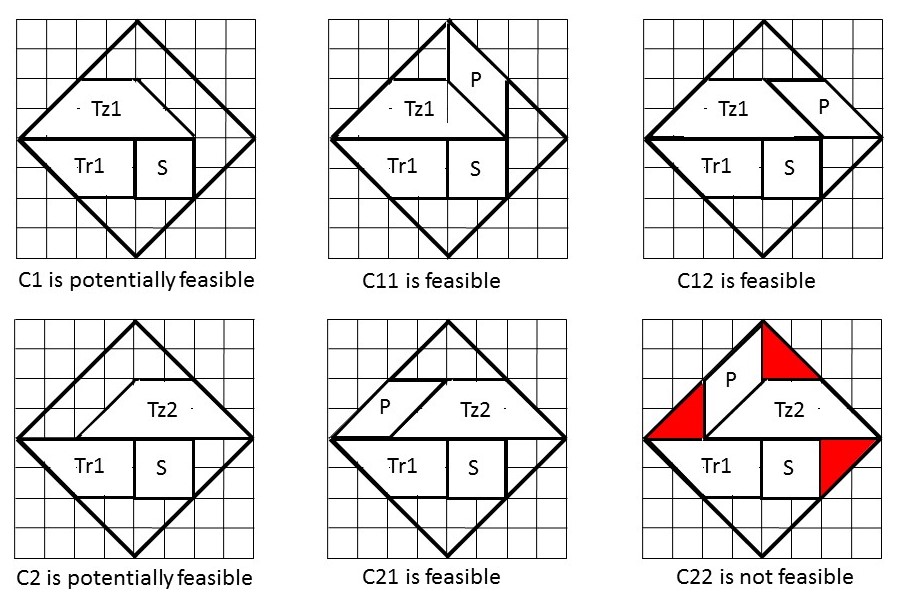}
\caption{All possible combinations when adding P to Tr1.Tz.S in $Square$}
\label{fig: Square_rem_feas}		%fig 14
\end{center}
\end{figure}
\clearpage
\begin{figure}[thb]
\begin{center}
\includegraphics[scale = 0.5]{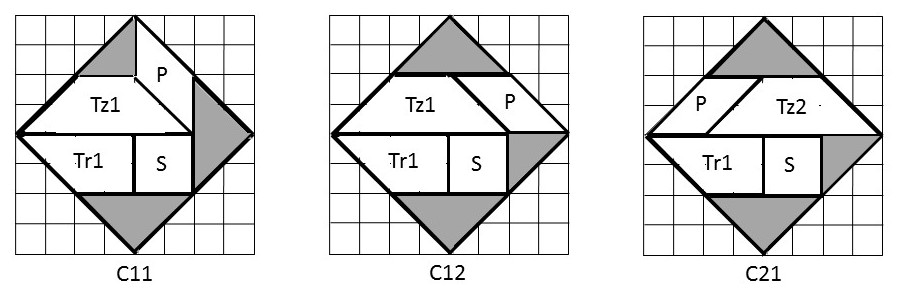}
\caption{The 3 feasible partitions for $Square$}
\label{fig: Square_all_feas}		%fig 15
\end{center}
\end{figure}

\noindent
Therefore, we return to the layouts C1 and C2 in Fig.~\ref{fig: Square_Tr1_Tz_S}. Let us try to add parallelogram P to them, see Fig.\ref{fig: Square_rem_feas}. Clearly, we have 2 options for adding P to $Square\_Tr1.Tz1.S$, see Figs.~\ref{fig: Square_rem_feas}-C11, C12. Clearly, after having added P we see that 3 empty regions are left in both C11 and C12, and they can be uniquely be covered by the triangles Tm1, Tm2 and Ts. Notice that the order of placing Tm1 and Tm2 is irrelevant.
Thus, layout $Square\_Tr1.Tz1.S$ can be completed in two different ways to a full feasible covering of $Square$.\newline
Similarly to C1 we can add P to $Square\_Tr1.Tz2.S$, but now only one layout (C21) is feasible, since in the other one (C22) we cannot add both triangles Tm1 and Tm2. See Figs.~\ref{fig: Square_rem_feas}-C21, C22. 
So, we have 
\begin{eqnarray}
	& & \hspace*{-0.7cm}\mbox{\rm\bf{Final Conclusion:}} 	
	\label{lab: Conclusion_3_layouts} \\		%concl 
	& & \hspace*{-0.7cm}\mbox{\rm\bf{Square $Tr1.Tz1.S$ can be fully covered by 3 different layouts, see Fig.~\ref{fig: Square_all_feas}. }} \nonumber
%		& & \hspace*{-0.7cm}\mbox{\rm\bf{see \; Fig.~\ref{fig: Square_all_feas}. }} \nonumber		
\end{eqnarray}
\clearpage

\section{The $Strip$ tangram $J14$}
\label{sec: Strip_J14}
We start by giving some preliminaries.
\subsection{Preliminaries}
\begin{figure}[thb]
\begin{center}
\includegraphics[scale = 0.6]{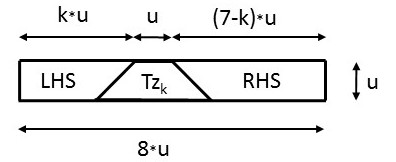}
\caption{Trapezium Tz with its LHS and RHS parts in strip $J14$.}
\label{fig: Strip_1}		%fig 16
\end{center}
\end{figure}
\begin{figure}[thb]
\begin{center}
\includegraphics[scale = 0.5]{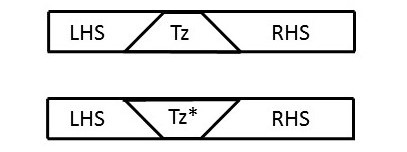}
\caption{Equivalence of the strips with Tz and Tz* (= Tz upside-down) in  $J14$}
\label{fig: Strip_2}		%fig 17
\end{center}
\end{figure}

\begin{figure}[thb]
\begin{center}
\includegraphics[scale = 0.55]{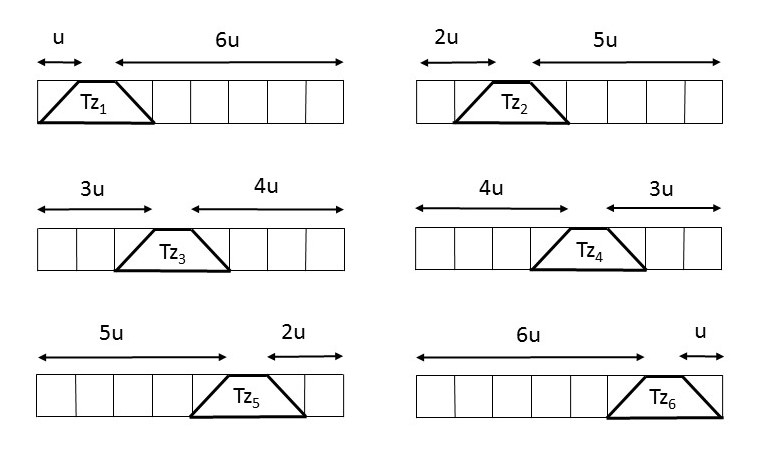}
\caption{All possible positions of $Tz$ in the strip.}
\label{fig: Strip_3}		%fig 18
\end{center}
\end{figure}

\begin{figure}[t!hb]
\begin{center}
\includegraphics[scale = 0.5]{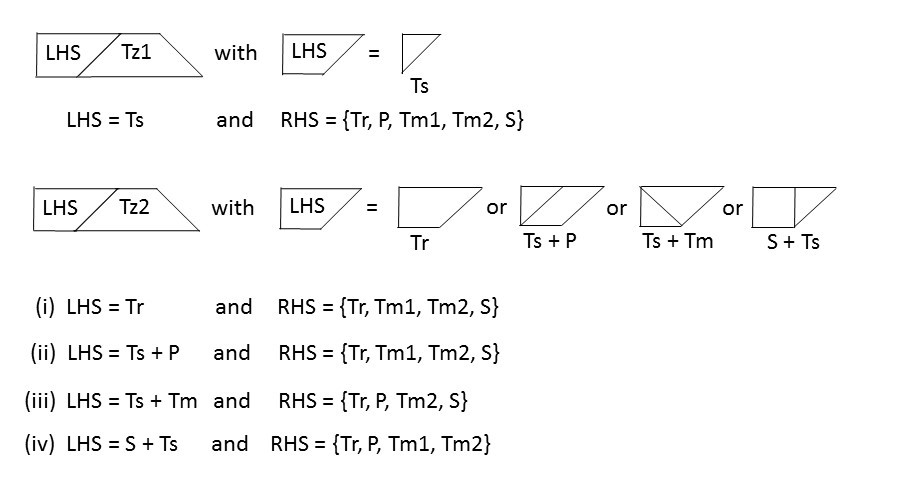}
\caption{All fillings for LHS and Tz1 and Tz2 in $Strip\_J14$.}
\label{fig: Jap14_LHS_RHS_Tz12}		%fig 19
\end{center}
\end{figure}

\begin{figure}[thb]
\begin{center}
\includegraphics[scale = 0.4]{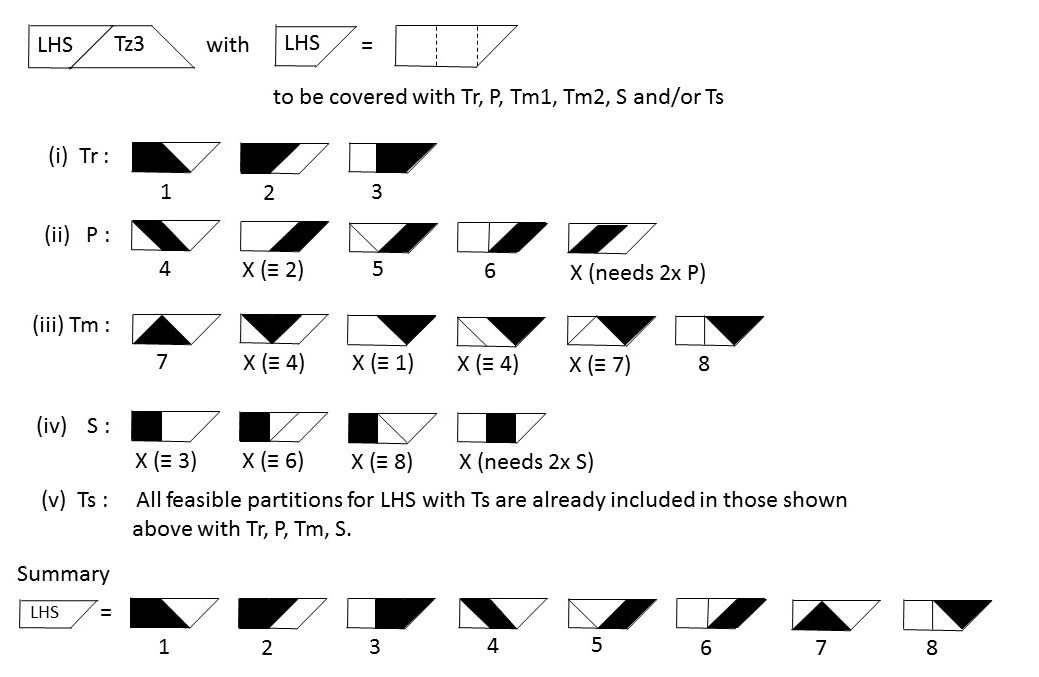}
\caption{All fillings for LHS and Tz3 in $Strip\_J14$.}
\label{fig: Jap14_LHS_RHS_Tz3}		%fig 20
\end{center}
\end{figure}

\noindent
Since Tz is an isosceles trapezium in a rectangular strip, we need to add tans at the left (LHS) as well as at the right hand side (RHS) of Tz. See Fig.\ref{fig: Strip_1}.
Clearly, to each side of Tz we have to add a $skew$-sided tan. Notice that the subscript $k$ in the name $Tz_k$ refers to the length of LHS.\newline
Furthermore, notice that each strip with Tz is equivalent to the strip with Tz upside-down (due to symmetry w.r.t. the bottom edge of the strip. See Fig.\ref{fig: Strip_2}.\newline\newline
Let us now consider all possible positions of Tz in the strip, see Fig.\ref{fig: Strip_3}. The strips will be identified by the Tz name. It is easily seen (by reflection w.r.t. the vertical left hand edge of the strip) that the following strips are equivalent: $Tz_4\, \equiv\, Tz_3,\; Tz_5\, \equiv\, Tz_2$ and $Tz_6\, \equiv\, Tz_1$. \newline
Thus, we need not to generate layouts with Tz4 up to Tz6. \newline
So, we have 
\begin{eqnarray}
	& & \hspace*{-0.7cm}\mbox{\rm\bf{Conclusion:}} 	
	\label{lab: Tz1_upto_Tz3} \\		%concl 
	& & \hspace*{-0.7cm}\mbox{\rm\bf{We only need to find all possible different layouts of the strip $J15$}} \nonumber \\
		& & \hspace*{-0.7cm}\mbox{\rm\bf{with $Tz_1$ up to $Tz_3$ \; (see Fig.~\ref{fig: Strip_3}). }} \nonumber		
\end{eqnarray}
Consider Fig.~\ref{fig: Strip_1} for $k=1,\;2,\;3$.
In Figs.~\ref{fig: Jap14_LHS_RHS_Tz12} and \ref{fig: Jap14_LHS_RHS_Tz3} we show all possible fillings for LHS corresponding to $Tz_k$. 
Clearly, once knowing all possible fillings for LHS and RHS we can find all different layouts of a partial strip $J14$ for each feasible combination of LHS and Tz (using Figs. 4 and 5) by adding step by step one tan from the set of remaining pieces. This can be done in a well-defined way by using the so-called backtracking procedure \cite{Tangram:ref_wikipedia_2017}, \cite{Tangram:ref_geeks_2017}. Below we will give more details. \newline 
\begin{figure}[thb]
\begin{center}
\includegraphics[scale = 0.35]{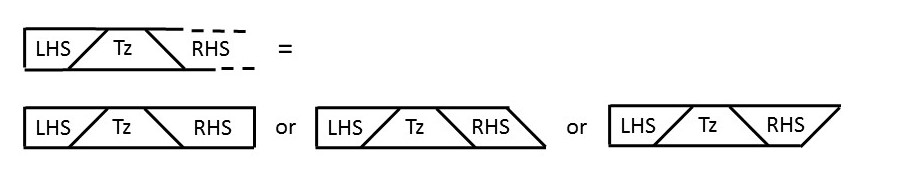}
\caption{The 3 possible rhs-edges of a partial $Strip\_J14$.}
\label{fig: Strip_Ends}		%fig 21
\end{center}
\end{figure} 
However, we first want to explain a few notational and visualization aspects.
\subsubsection{Visualisation aspects}
\begin{figure}[thb]
\begin{center}
\includegraphics[scale = 0.52]{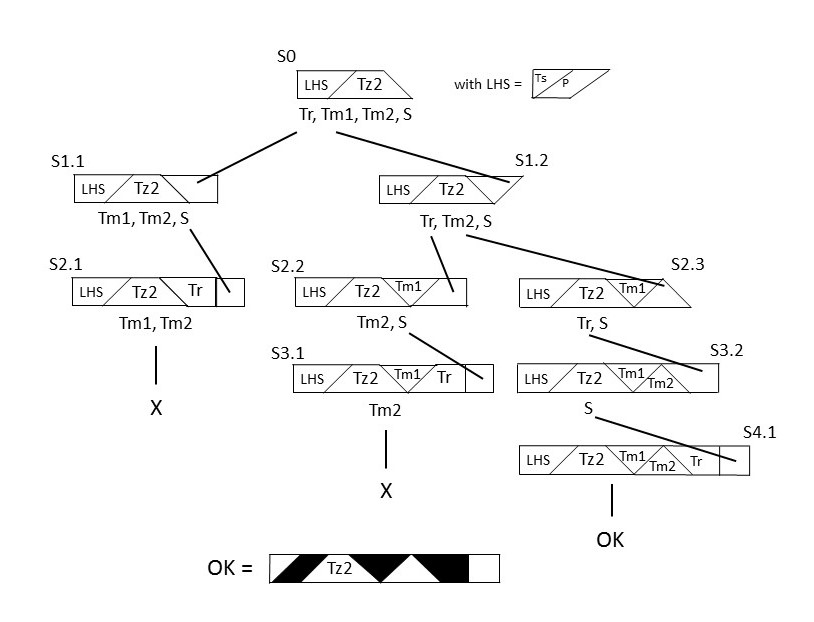}
\caption{Example of a tree structure for finding all layouts of the strip J14.}
\label{fig: Tree_Visual1}		%fig 22
\end{center}
\end{figure}

We start with a partial strip LHS+Tz and add a suitable tan (called T1 for simplicity) next to Tz. We get the (partial) strip $LHS+Tz+RHS = LHS+Tz+T1$. Notice that there are 2 options for each tan to be added: either the tan fits to the strip in a unique way or the tan does not fit at all. \newline
Clearly, the rhs-edge of each partial strip is vertical, left- or right skew, see Fig.~\ref{fig: Strip_Ends} and the final strip must have a vertical rhs-edge. By adding more tans $T2,\; T3\; , \cdots $ the size of the strip grows. This process can be visualized by a tree structure. It is easily seen that the tree will grow widely when showing all partial strips in full detail. However, we will show not all details but only the actual  strip globally. \newline
We will clarify this by the following representative example. \newline\newline
$\mathbf{Visualisation\; of\; the\; backtracking\; process}$\newline\newline
Consider Fig.~\ref{fig: Tree_Visual1} where we have a partial strip S0 = LHS+Tz with LHS = Ts+P. \newline
Then the remaining tans are Tr, Tm1, Tm2 and S and these names are listed under S0. We can extend S0 with either Tr or Tm, resulting in the partial strips S1.1 and S1.2, respectively. \newline
Clearly, S1.1 is fully rectangular while S1.2 has a skew rhs-edge. The remaining tans for extending S1.1 are Tm1, Tm2 and S (these names are listed under S1.1). However, S1.1 can only be extended with S, resulting in S2.1, with Tm1 and Tm2 left. It is easily seen that these tans cannot be used anymore for extending  S2.1. This fact is indicated by X in the figure and this branch of the tree ends here. So, adding the tans in this order does not result in a full strip. \newline
On the other hand, when first adding one of the triangles Tm1, Tm2 (say Tm1) to S0 we find the strip S1.2, which can be extended by adding either Tr or Tm2, resulting in the strips S2.2 and S2.3, respectively. Now the remaining tans for extending S2.2 are Tm2 and S (see below S2.2), but only S can be added, resulting in S3.1. However, S3.1 cannot be extended further, again indicated by X and this branch also ends. \newline
Let us now consider S2.3. Here the remaining tans for extension are Tr and S. Clearly, S2.3 can be extended with only Tr, giving strip S3.2 and tan S is left. Finally, S3.2 can be extended by S, giving a final full strip, indicated by OK. The structure of the final full strip is also given. \newline\newline
\begin{figure}[thb]
\begin{center}
\vspace{2.cm}
\hspace{-2.4cm}
\includegraphics[scale = 0.55]{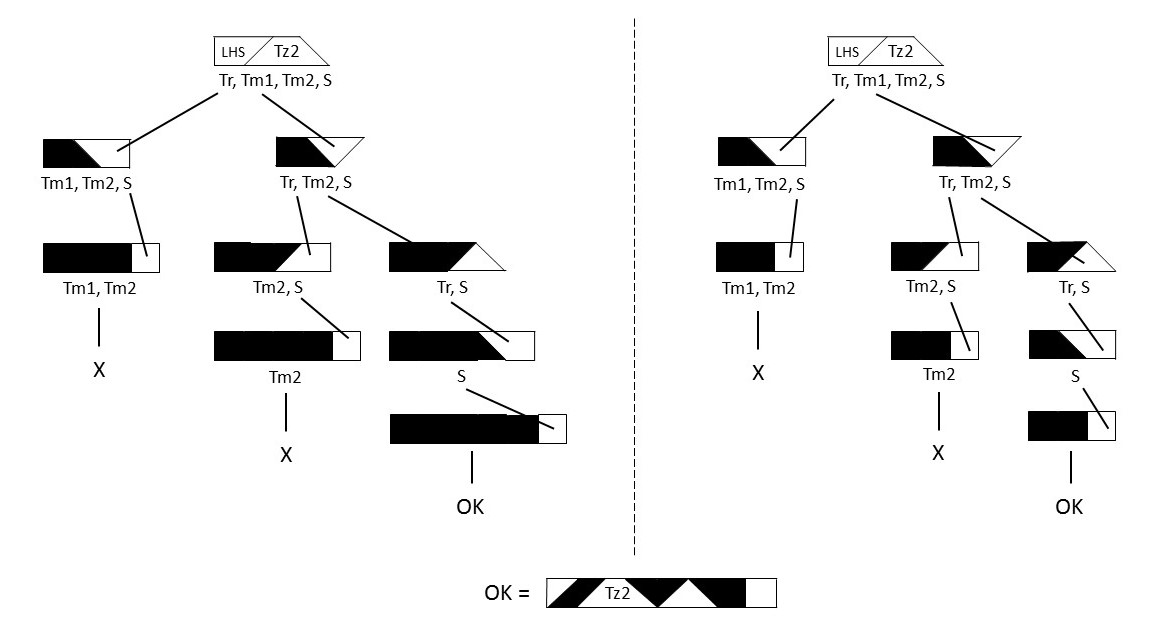}
\caption{Two global visualizations of a tree structure for J14.}
\label{fig: Tree_Visual2}		%fig 23
\end{center}
\end{figure}
$\mathbf{Simplification\; of\; the\; tree\; visualization}$
\newline \newline
It should be noticed that in fact we do not need to know the full filling details of the intermediate partial strips. Indeed, only the overall shape of each partial strip is relevant for finding a possible extension. So, we can replace the detailed scheme above by a more global scheme, as shown in Fig.~\ref{fig: Tree_Visual2}-left. We emphasize that in this tree at each level (except top level) we have only indicated (in black) the $global$ shape (and actual) size of the partial strip in the previous level. However, the newly added tan at the current level is shown with its $actual$ shape and $actual$ size.\newline\newline
Finally, we even can reduce the width of the figure a bit more by using the $same$ size for all partial strips, but preserving their global shape, as shown in Fig.~\ref{fig: Tree_Visual2}-right. \newline
In all next figures we will only use the latter compact visualization.
\newline\newline
\subsection{Finding all possible layouts of strip $J14$}
\subsubsection{The layouts of the strip with Tz1 and Tz2}
\begin{figure}[thb]
\begin{center}
\hspace{-2cm}
\includegraphics[scale = 0.42]{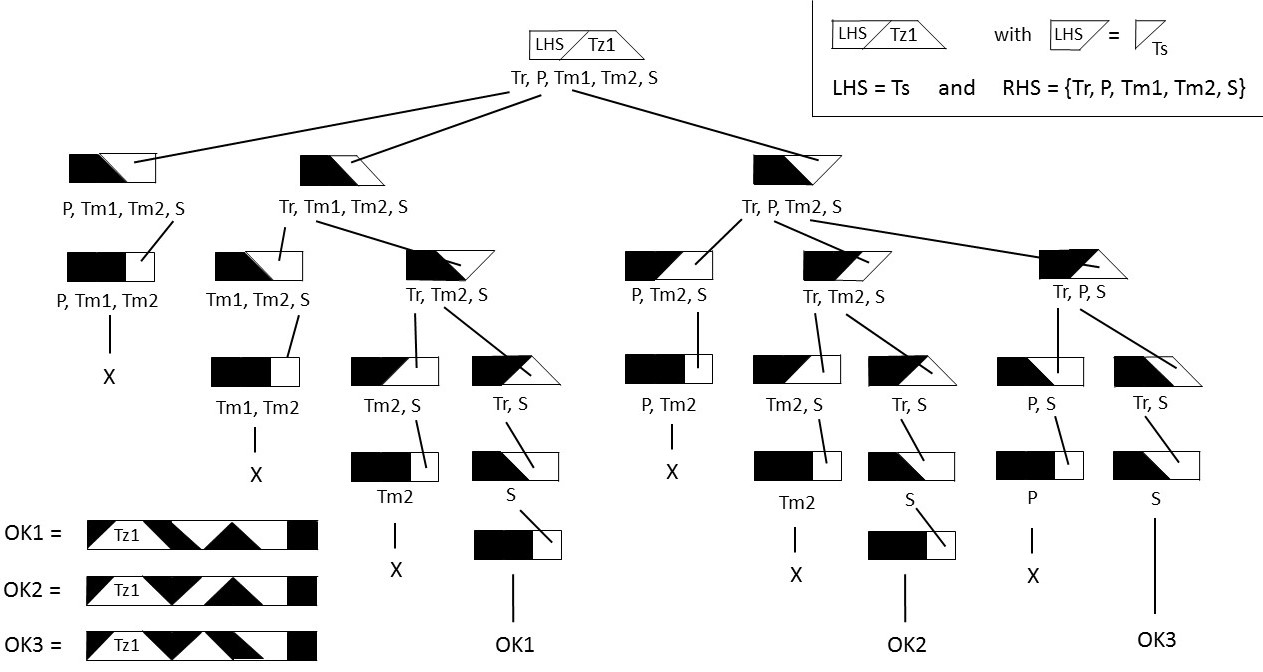}
\caption{Visualization of Strip $J14\_Tz1$.}
\label{fig: J14_Tz_1}		%fig 24
\end{center}
\end{figure}

As indicated in Fig.~\ref{fig: Jap14_LHS_RHS_Tz12} we see that (i) in case of Tz1 the corresponding LHS consists of one single tan, being Ts, and  (ii) in case of Tz2 that LHS consists of one or two tans. All layouts for $J14$ with $Tz1$ and $Tz2$ are given in Figs.~\ref{fig: J14_Tz_1} up to \ref{fig: J14_Tz_24}.
\clearpage
\subsubsection{The layouts of $J14$ with $Tz1$}
\subsubsection{The layouts of $J14$ with $Tz2$}
\begin{figure}[thb!]
\begin{center}
\hspace{-2.5cm}
\includegraphics[scale = 0.42]{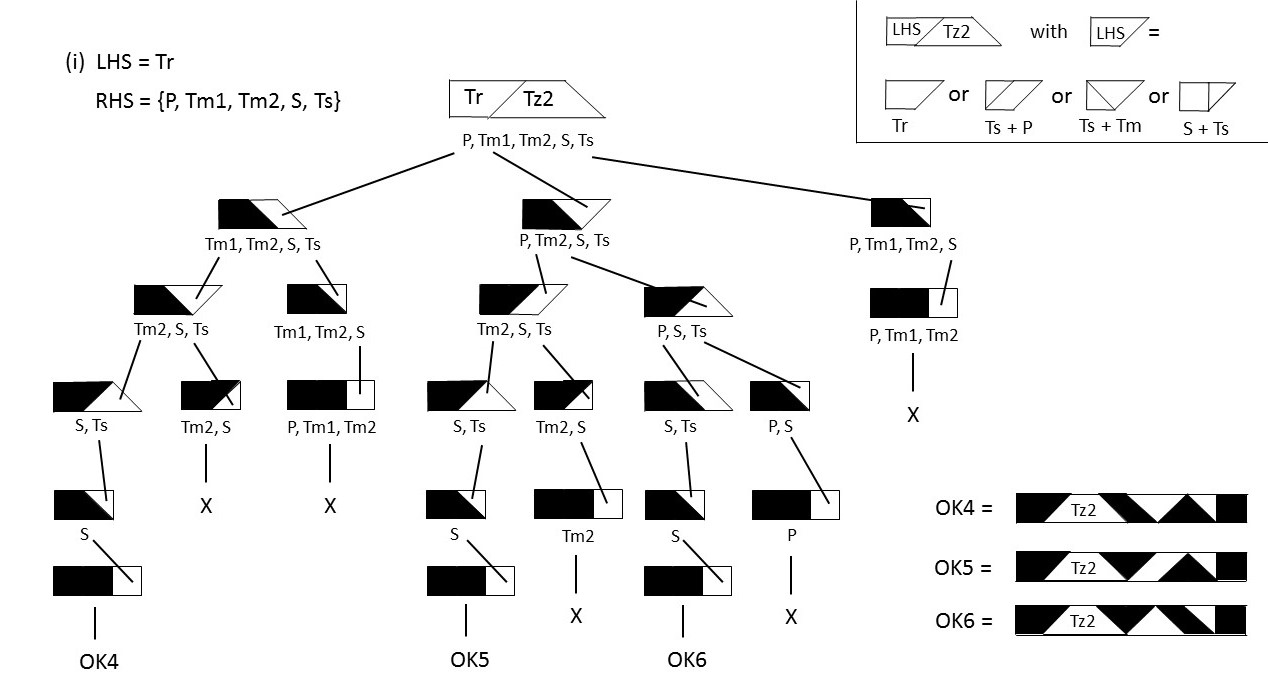}
\caption{Visualization of Strip $J14\_Tz2.1$.}
\label{fig: J14_Tz_21}		%fig 25
\end{center}
\end{figure}
\clearpage
\begin{figure}[thb]
\begin{center}
\hspace{-2cm}
\includegraphics[scale = 0.5]{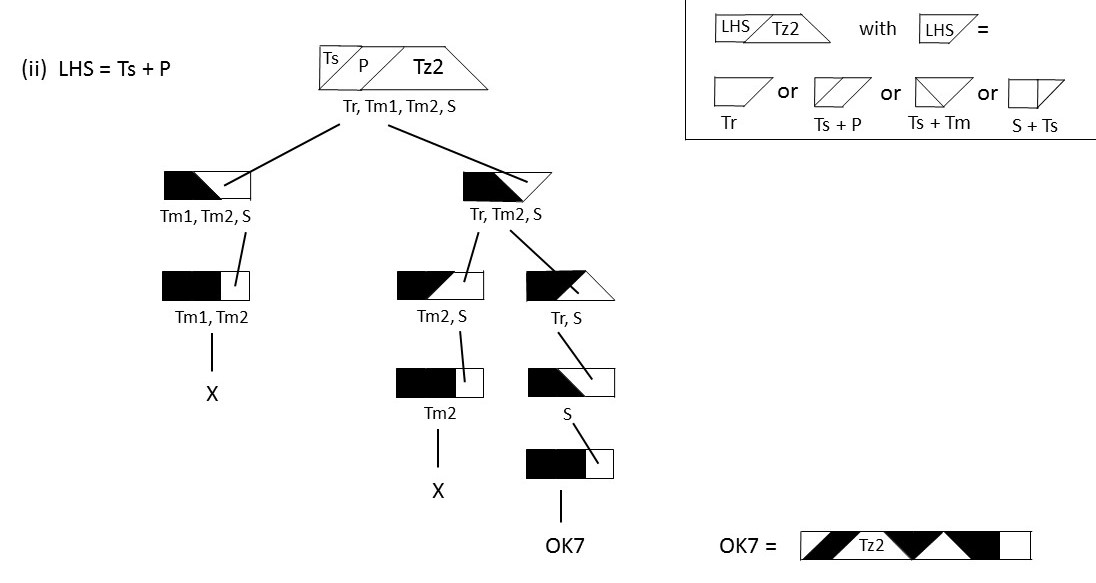}
\caption{Visualization of Strip $J14\_Tz2.2$.}
\label{fig: J14_Tz_22}		%fig 26
\end{center}
\end{figure}
\begin{figure}[thb!]
\begin{center}
\hspace{-2cm}
\includegraphics[scale = 0.5]{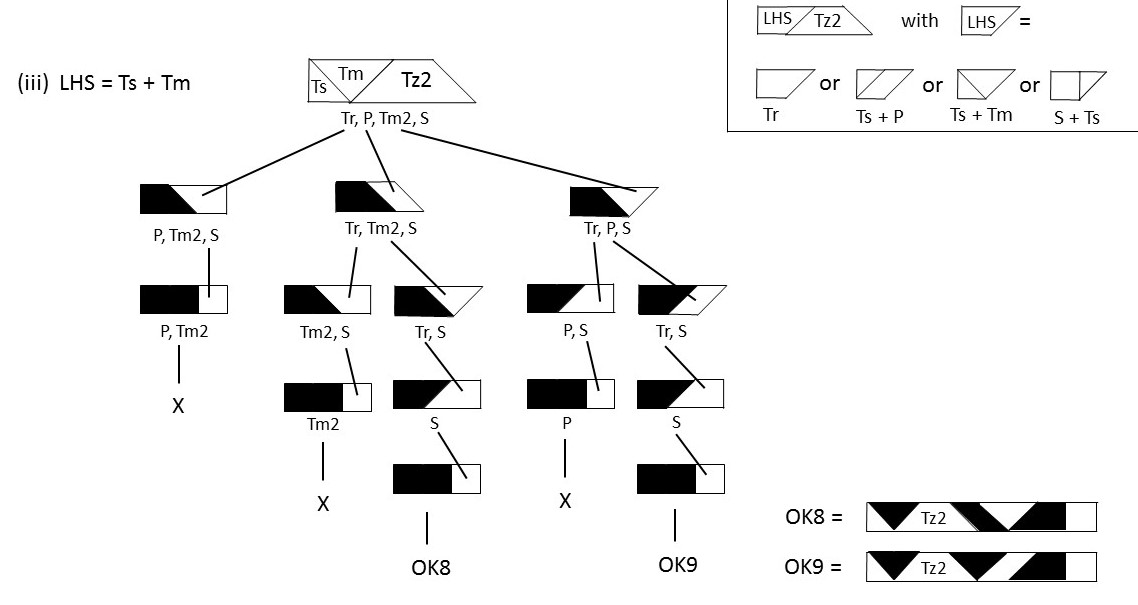}
\caption{Visualization of Strip $J14\_Tz2.3$.}
\label{fig: J14_Tz_23}		%fig 27
\end{center}
\end{figure}
\clearpage
\begin{figure}[thb]
\begin{center}
\hspace{-2cm}
\includegraphics[scale = 0.45]{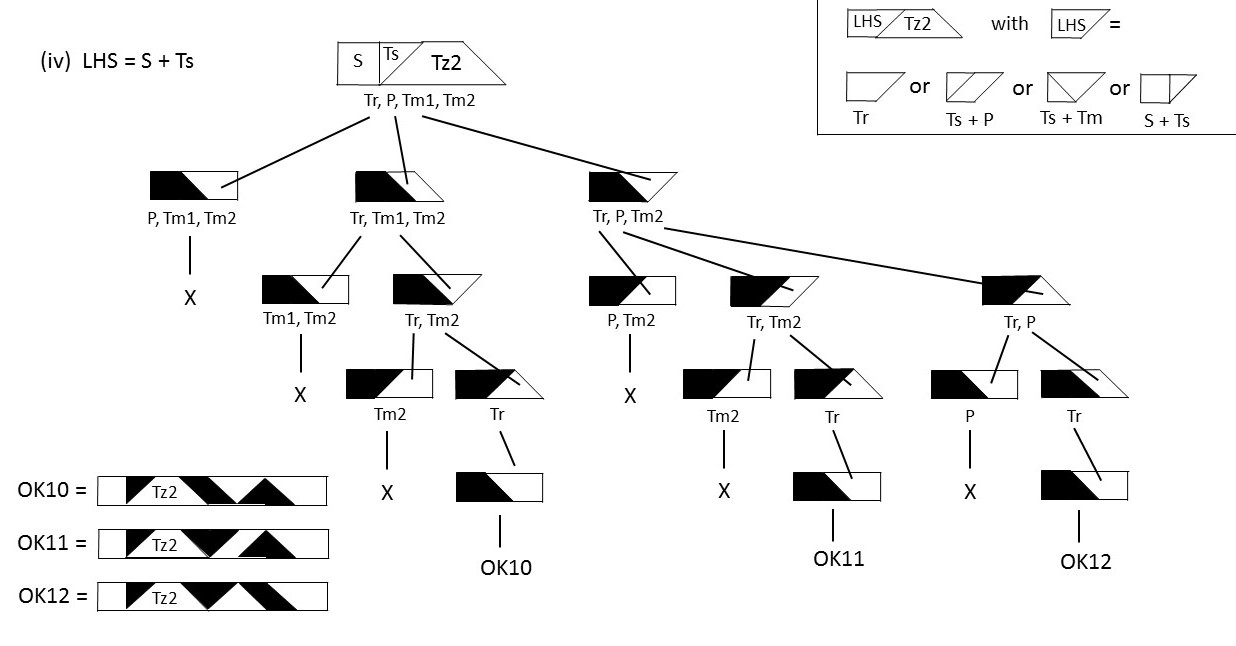}
\caption{Visualization of Strip $J14\_Tz2.4$.}
\label{fig: J14_Tz_24}		%fig 28
\end{center}
\end{figure}
\subsubsection{The layouts of $J14$ with $Tz3$}
\begin{figure}[thb!]
\begin{center}
\hspace{-2cm}
\includegraphics[scale = 0.44]{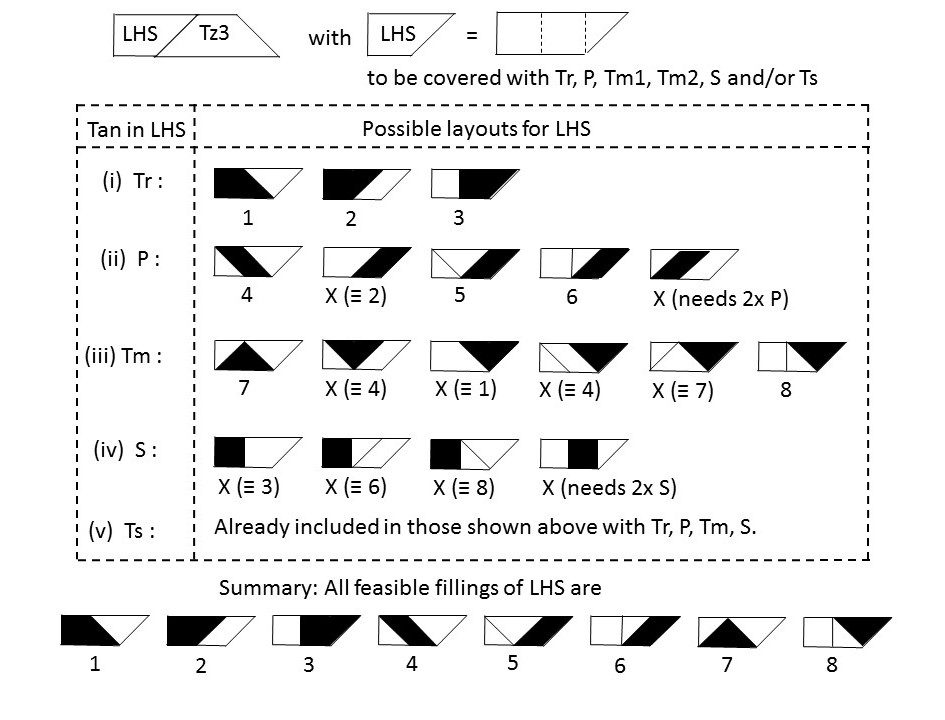}
\caption{Investigation of the LHS partitions of Strip $J14\_Tz3$.}
\label{fig: J14_Tz_29}			%fig 29  = table
\end{center}
\end{figure}
\clearpage
\subsubsection{The layouts of the strip with Tz3}
As indicated in Fig.~\ref{fig: J14_Tz_29} we see that in case of Tz3 the corresponding LHS consists of 2 or 3 tans, resulting in 8 feasible fillings. Notice that the left-edge of LHS is vertical, to be realized by one of the 3 tans Tr, S and Ts. Further, RHS has a skew left edge and a vertical right edge. The latter must also be realized by one of the 3 tans Tr, S and Ts. So, when fixing one of these 3 for realizing the vertical edge of LHS, at most 2 of them are left for building RHS.  
\begin{figure}[thb]
\begin{center}
\hspace{-2.5cm}
\includegraphics[scale = 0.45]{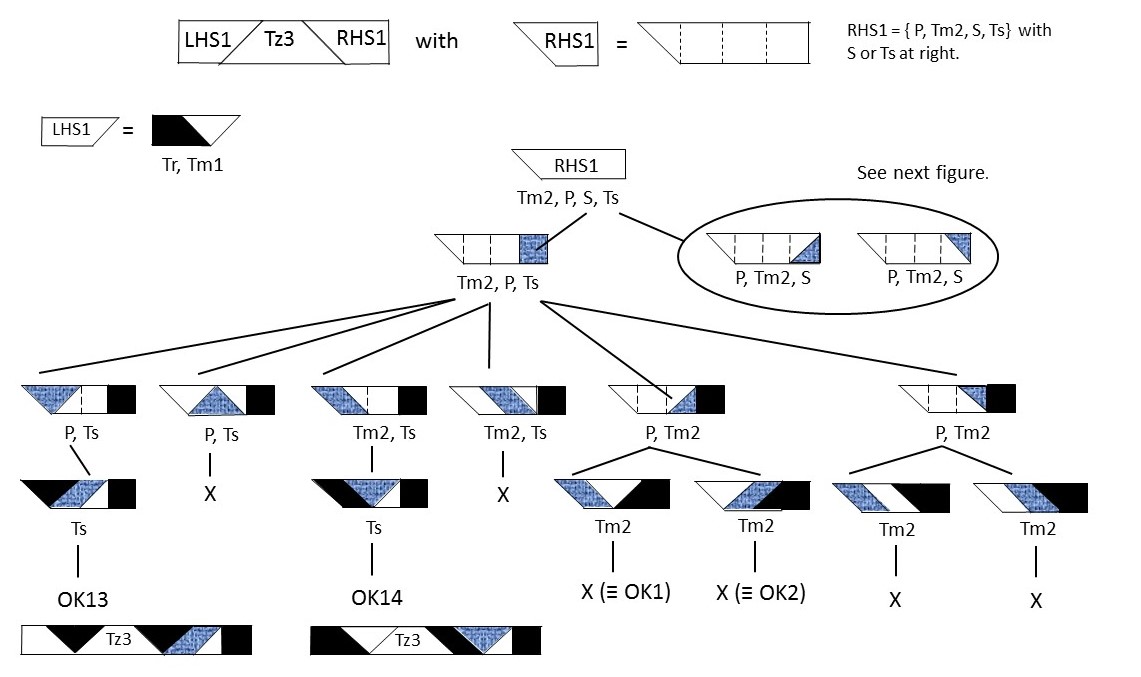}
\caption{Investigation of Strip $J14\_Tz31$, part (i).}
\label{fig: J14_Tz_31a}			%fig 30
\end{center}
\end{figure}
\begin{figure}[thb]
\begin{center}
\vspace{1cm}
\hspace{-2cm}
\includegraphics[scale = 0.48]{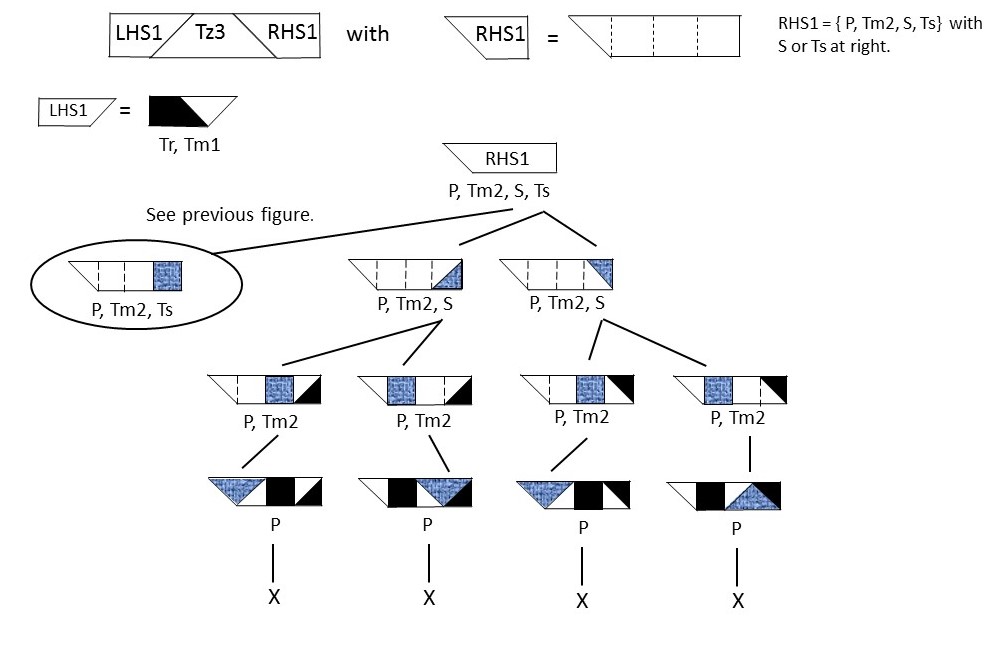}
\caption{Investigation of Strip $J14\_Tz31$, part (ii).}
\label{fig: J14_Tz_31b}			%fig 31
\end{center}
\end{figure}
\newline\newline
In Figs.~\ref{fig: J14_Tz_31a} and \ref{fig: J14_Tz_31b} we show the partial fillings of RHS ending at right with S or Ts.\newline
Clearly, Fig.~\ref{fig: J14_Tz_31a} is self-explaining.\newline
Fig.~\ref{fig: J14_Tz_31b} shows the fillings of RHS ending at right with Ts. Next, after having added S in the 4 RHS-strips we still need to add P and Tm. However, for a full filling of all these 4 sub-strips we need an additional tan Ts, but this is not available. Thus, this branch of the tree does not end with a feasible filling of RHS.
\newline
Figs.~\ref{fig: J14_Tz_32} and \ref{fig: J14_Tz_33} are also self-explaining.
\begin{figure}[thb]
\begin{center}
\hspace{-2cm}
\includegraphics[scale = 0.5]{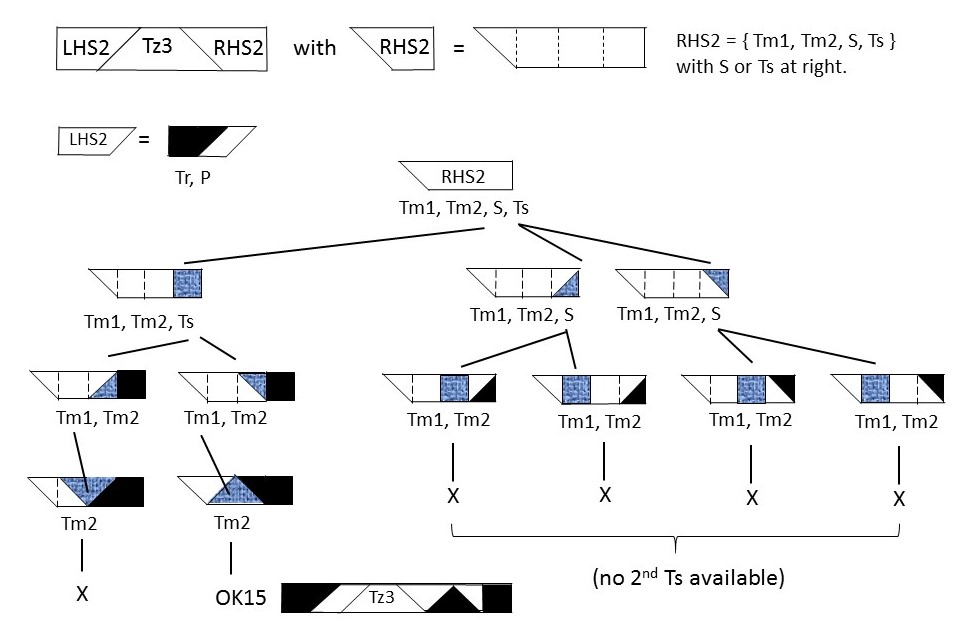}
\caption{Investigation of Strip $J14\_Tz32$.}
\label{fig: J14_Tz_32}		%fig 32
\end{center}
\end{figure}
\begin{figure}[thb]
\begin{center}
\hspace{-2cm}
\includegraphics[scale = 0.5]{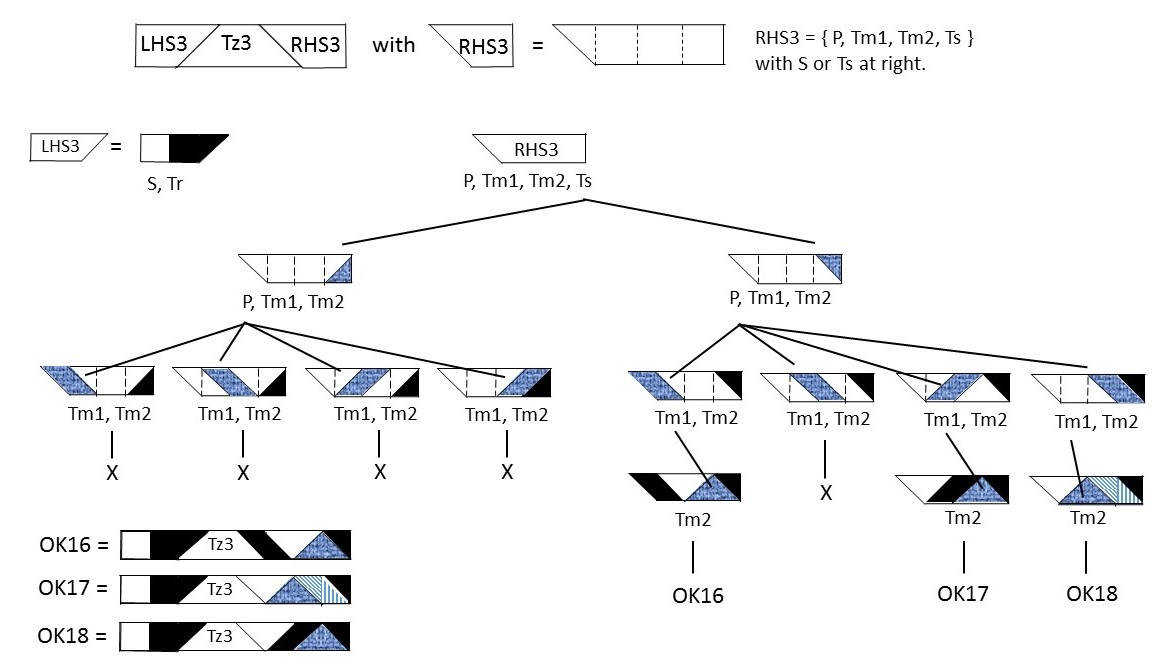}
\caption{Investigation of Strip $J14\_Tz33$.}
\label{fig: J14_Tz_33}		%fig 33
\end{center}
\end{figure}
\newline\newline
In Fig.~\ref{fig: J14_Tz_34} we show the partial fillings of RHS4 ending at right with Tr or S.
\begin{figure}[thb]
\begin{center}
\hspace{0 cm}
\includegraphics[scale = 0.6]{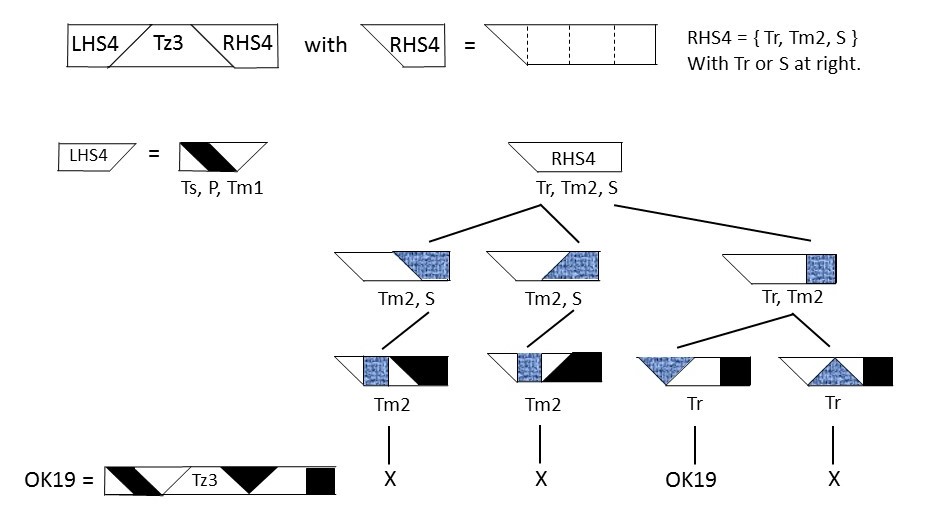}
\caption{Investigation of Strip $J14\_Tz3.4$.}
\label{fig: J14_Tz_34}			%fig 34
\end{center}
\end{figure}
%
%\clearpage
In Fig.~\ref{fig: J14_Tz_35} we consider the case where LHS5 consists of the same tans (Ts, P and Tm1) as in the previous figure.  Consequently, RHS5 is identical to RHS4 since the same tans for RHS5 are available as for RHS4. However, since LHS5 has a partition being different from that of LHS4, we find with LHS5, Tz3 and RHS5 a different partition of the whole strip which is different from that in the previous case with LHS4, Tz3 and RHS4.
\begin{figure}[thb]
\begin{center}
\hspace{-2cm}
\includegraphics[scale = 0.6]{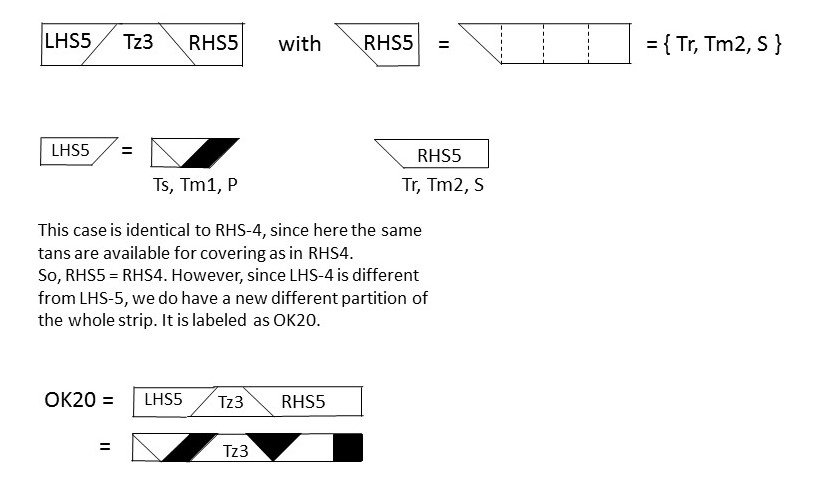}
\caption{Investigation of Strip $J14\_Tz3.5$.}
\label{fig: J14_Tz_35}			%fig 35
\end{center}
\end{figure}
\newline\newline
In Fig.~\ref{fig: J14_Tz_36} we consider the case LHS6 with Tr, Tm1 and Tm2. Since both triangles do not have a vertical edge we have to place Tr at right, resulting in 2 layouts. Clearly, only one feasible partition can be made with Tm1 and Tm2. 
\clearpage
\begin{figure}[thb]
\begin{center}
\hspace{-2cm}
\includegraphics[scale = 0.58]{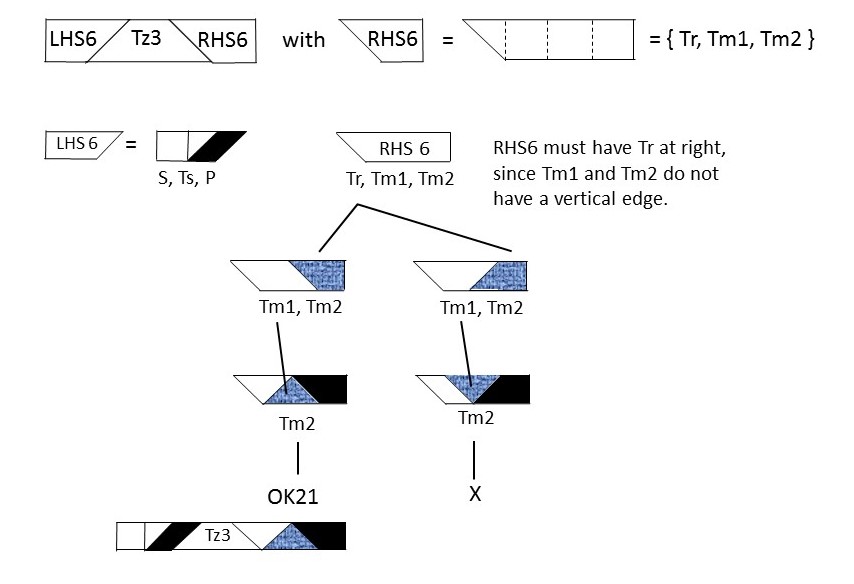}
\caption{Investigation of the Strip $J14\_Tz3.6$.}
\label{fig: J14_Tz_36}			%fig 36
\end{center}
\end{figure}
\begin{figure}[thb!]
\begin{center}
\hspace{0 cm}
\includegraphics[scale = 0.58]{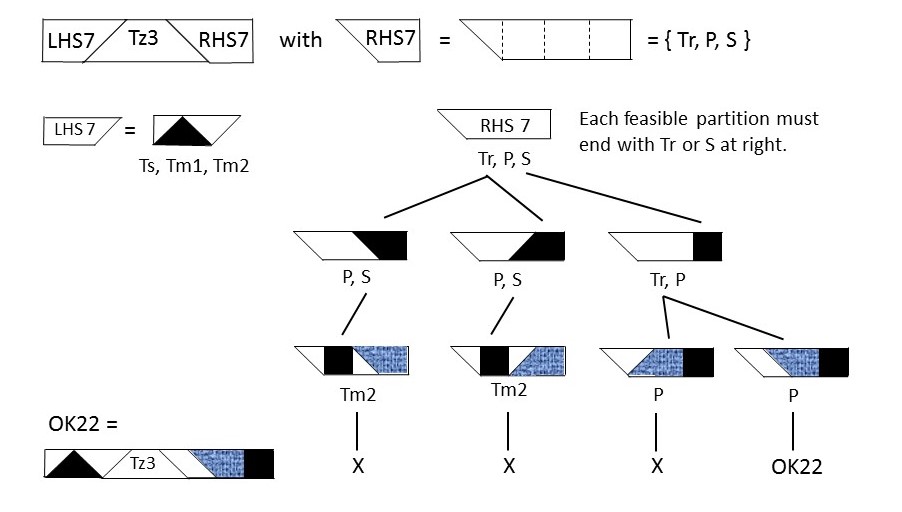}
\caption{Investigation of the Strip $J14\_Tz3.7$.}
\label{fig: J14_Tz_37}		%fig 37
\end{center}
\end{figure}

\begin{figure}[thb]
\begin{center}
\hspace{0 cm}
\includegraphics[scale = 0.45]{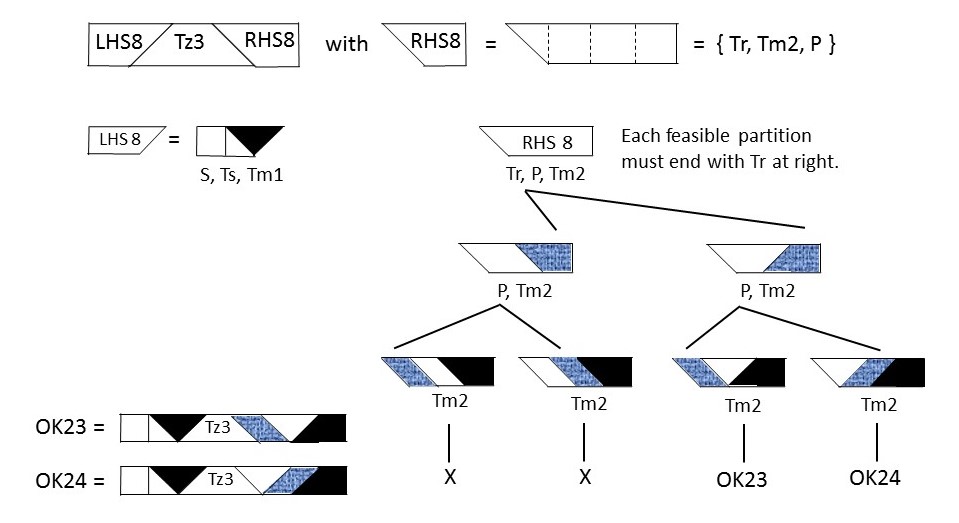}
\caption{Investigation of te Strip $J14\_Tz3.8$.}
\label{fig: J14_Tz_38}			%fig 38
\end{center}
\end{figure}

\noindent
In Fig.~\ref{fig: J14_Tz_37} we consider case LHS7 consisting of Ts, Tm1 and Tm2. Hence, RHS7 must contain Tr, P and S. Both Tr and S have a vertical edge, so these tans must be placed at right, resulting in 3 partial fillings for RHS. Next  we can add S and Tr to these layouts. Notice that after having added S, two tans of type Ts are required for a full filling, but these are not available. Similarly, after having added Tr we find one layout where only Tm can be added while only P is present. The other layout we have to add P and tin this case this is possible, resulting in a feasible layout (denoted by OK22).  
\newline\newline
Next we have to study the strip with LHS8 with S, Ts and Tm1. Then RHS8 must contain Tr, P and Tm2. Clearly, Tr must be placed at right of RHS8 since P and Tm2 do not have a vertical edge. This gives two options where we can add P, resulting in 4 partial strips. Finally, herein Tm has to be included. This is only possible in two strips. Hence, this case with LHS8 and RHS8 results in 2 feasible partitions for the whole  strip. See Fig.~\ref{fig: J14_Tz_38}.
\subsubsection{Summary of the analysis of $Strip\_ J14$}
\begin{figure}[t!hb]
\begin{center}
\hspace{0cm}
\includegraphics[scale = 0.65]{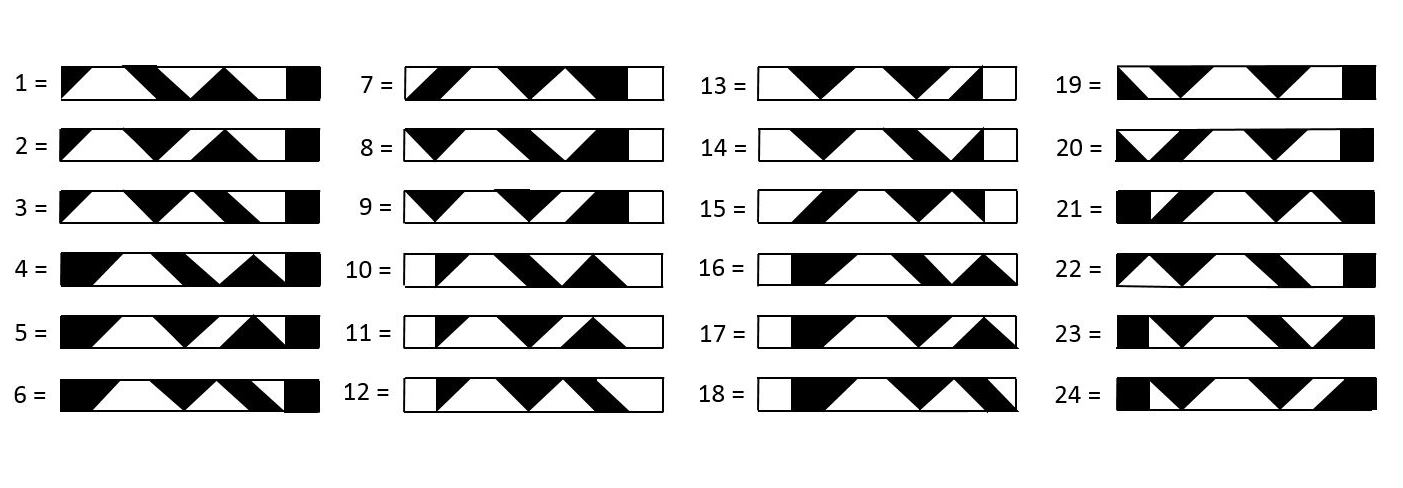}
\caption{Strip $J14$ with all its 24 different partitions (``handmade'').}
\label{fig: J14_all_solutions}			%fig 39
\end{center}
\end{figure}

\begin{figure}[t!hb]
\begin{center}
%\hspace{-0.8cm}
%\includegraphics[scale = 0.6, angle=90]{J14_all_solutions}
\includegraphics[scale = 0.45]{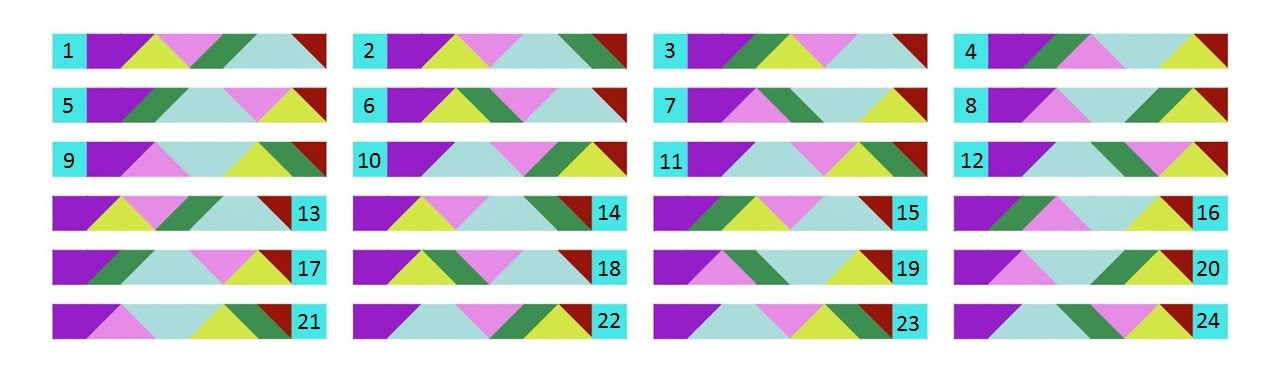}
\caption{Strip $J14$ with the 24 partitions (``computer generated'').}
\label{fig: japsolred_14_nrs}			%fig 40
\end{center}
\end{figure}

The investigations above for finding all possible partitions of $Strip J14$ with the 7 given tans can be summarized by the following steps:
\begin{enumerate}
\item Determine all feasible positions and orientations of the isosceles trapezium Tz in strip $J14$, see Fig.~\ref{fig: Strip_3};
\item Determine all possible partitions for the LHS in the strip, see Figs.~\ref{fig: Jap14_LHS_RHS_Tz12}, \ref{fig: Jap14_LHS_RHS_Tz3} and \ref{fig: J14_Tz_29};
\item For each of the LHS-partitions: determine all possible layouts for the right hand side (RHS) in the strip using the backtracking principle and being visualized by a tree stucture, see Figs.~\ref{fig: Tree_Visual1} up to \ref{fig: J14_Tz_38}.
\end{enumerate} 
In this way we find all possible partitions of the complete strip $J14$, giving in total 24 different solutions. They are shown in Fig.~\ref{fig: J14_all_solutions}. Moreover, these 24 solutions have also been found by our computer program (see Fig.~\ref{fig: japsolred_14_nrs} and section \ref{sec: backtr_algo}). \newline
For convenience, the equivalence between the partitions (``handmade'' and ``computer generated'') in both figures is given in Table \ref{Correspondences_J14_sols} after Fig.~\ref{fig: japsolred_14_nrs}.
Notice that sometimes we have to apply a horizontal or vertical reflection to a particular partition in Fig.~\ref{fig: J14_all_solutions} to find the same picture in Fig.~\ref{fig: japsolred_14_nrs}.
%
%%			add correspondence table here
\begin{table}[t!hb]
\caption{Correspondence between the handmade and computer generated solutions for Strip $J14$.}
\vspace{2mm}
\centering
\begin{tabular}{|c |c c c c c c c c c c c c |}
\hline
%	  		   &  &   &   &    &    &    & 		&   &    &    &    & 		\\	%empty line
	$Handmade\; (part 1)$ 
	&$1$  &$2$  &$3$  &$4$   &$5$   &$6$ &$7$ &$8$
						   &$9$  &$10$ &$11$ &$12$ \\
%	  		   &  &   &   &    &    &    & 		&   &    &    &    & 		\\	%empty line
 $Computer\; generated$ 
 &1 &6 &3 &24 &22 &23 &2 &7 &4 &13 &18 &15 \\
\hline
%	  		   &  &   &   &    &    &    & 		&   &    &    &    &  \\		
%empty line	  		   
	$Handmade\; (part 2)$ 
	&$13$  &$14$  &$15$   &$16$ &$17$ &$18$
						   &$19$  &$20$ &$21$ &$22$ &$23$  &$24$\\
%	    &  &   &   &  &   &    & 	&   &    &    &    & 		\\	
 $Computer\; generated$ 
 &21 &20 &17 &12 &10 &11 &9 &8 &14 &5 &19 &16 \\
\hline
\end{tabular}
\label{Correspondences_J14_sols}
\end{table}
\clearpage
\subsubsection{Strip $J14$ with its twin layouts}
\begin{figure}[thb!]
\begin{center}
\hspace{-2cm}
\includegraphics[scale = 0.55]{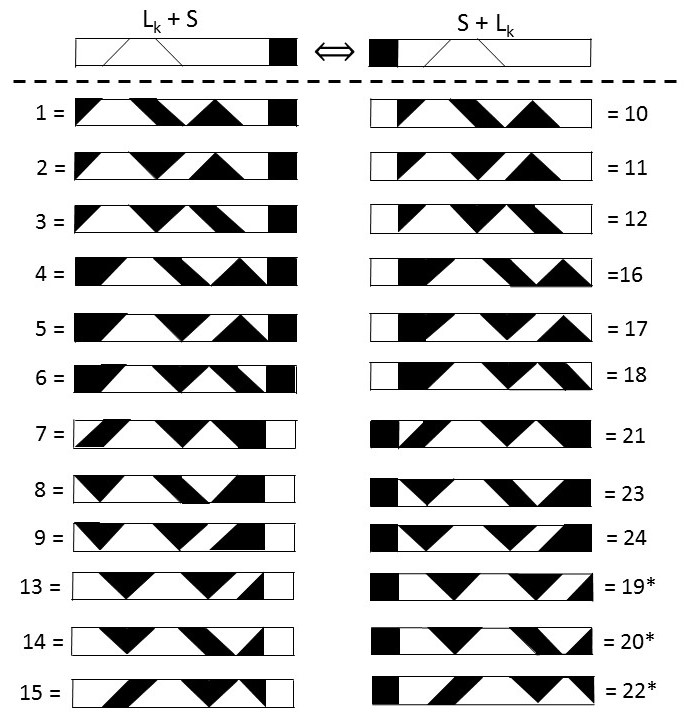}
\caption{Strip $J14$ with its 12 twin layouts. \newline The strips 19*, 20* and 22* are the horizontally flipped strips 19, 20, and 21 in Fig.~\ref{fig: J14_all_solutions}.}
\label{fig: J14_twin_solutions}			%fig 41
\end{center}
\end{figure}

Let us consider the 24 solutions in Fig.\ref{fig: J14_all_solutions} in more detail. We can divide each of the strips into the tan S and a 7S-wide substrip L. Notice that S is either at the left or at the right side of the full strip. We will denote the 24 full strips by $F_1$ up to $F_{24}$ and their substrips by $L_1$ up to $L_{24}$. 
Thus, we have either $F_k = L_k +S$ or $F_k = S + L_k$ for 
$k=1, \cdots , 24$. 
The strips $L_k +S$ and $S + L_k$ will be called twins and their twin-relationship will be indicated by $L_k + S \Leftrightarrow S + L_k$. 
It is easily seen from Fig.~\ref{fig: J14_all_solutions} that we have the following twins, shown in Fig.~\ref{fig: J14_twin_solutions}. 
\clearpage
\subsection{Finding all different partitions for $J15$ and $J16$ by a combinatorial approach}
\begin{figure}[thb]
\begin{center}
\includegraphics[scale = 0.5]{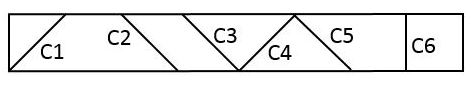}
\caption{The cuttings edges of strip $J14$-$1$.}
\label{fig: J14_cuttings}			%fig 42
\end{center}
\end{figure}

We start with considering the structure of all feasible partitions found for the 12 twin pairs in Fig.~\ref{fig: J14_twin_solutions}.
Notice that inside each of the 24 partitions there are 6 joint edges for each pair of adjacent tans. In particular, precisely 5 of these joint edges are skew, and only one joint (with S) is vertical. We can cut $J14$ along each of these cutting edges, resulting in 2 separate sub-strips. Next we can reverse the order of these sub-strings and glue them together. Apparently, the latter can be done in two ways (i.e., in original or in upside-down orientation) when the cutting edge is skew.\newline
We will discuss the possible situations for all strips. This is done by using one representative strip. To this end, we can take the first strip $J14$-$1$ in Fig.\ref{fig: J14_all_solutions}. \newline
The cutting edges will be denoted by $C1$ up to $C6$, with $C6$ being vertical. See Fig.~\ref{fig: J14_cuttings}.
\begin{figure}[thb]
\begin{center}
%\hspace{-3cm}
\includegraphics[scale = 0.55]{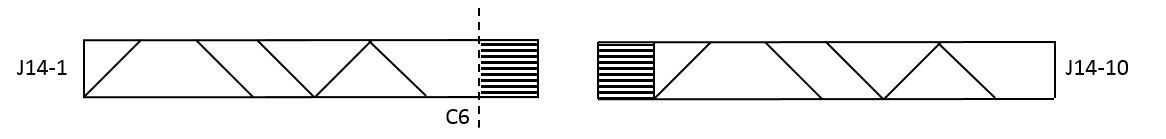}
\caption{Vertical cutting and pasting of strip $J14$-$1$ gives strip $J14$-$10$ .}
\label{fig: J14_cuttings2}			%fig 43
\end{center}
\end{figure}
\begin{figure}[thb]
\begin{center}
\hspace{-1.cm}
\includegraphics[scale = 0.48]{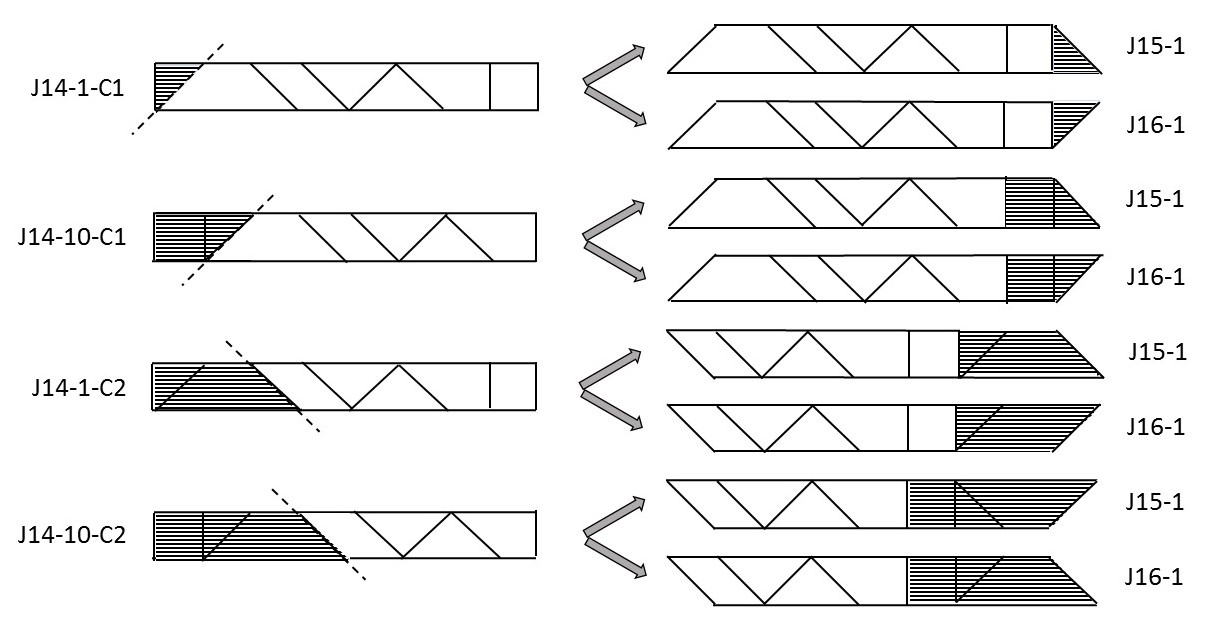}
\caption{Cutting and pasting the twin strips $J14$.}
\label{fig: J14_twin_cuttings}			%fig 44
\end{center}
\end{figure}
\newline\newline
Let us first consider the case with a vertical cutting edge (i.e. the cutting along vertical edge $C6$ of tan S inside the strip $J14$-$1$. This is illustrated in Fig.~\ref{fig: J14_cuttings2}-Left. 
We obtain the tan S and a sub-strip of width 7 S. After glueing S in front of the sub-string we get strip $J14$-${10}$, see Fig.~\ref{fig: J14_cuttings2}-Right. 
It is easily seen that carrying out this process for all strips in the lhs column of Fig.~\ref{fig: J14_twin_solutions} results in the creation of all corresponding twin strips in its rhs column.
\newline\newline
Next we consider the case with a skew cutting edge ($C1$ up to $C5$). The process of ``cut and paste''  (in two ways) is visualized by the two (representative) examples $J14$-$C1$ and $J14$-$C2$ in Fig.~\ref{fig: J14_twin_cuttings}.
Then we obtain two new strips, one having the shape of $J15$, and the other one that of $J16$. 
\newline\newline
Recalling Fig.~\ref{fig: J14_twin_solutions} we know that $J14$-$1$ and $J14$-$10$ are twins. We can apply the cut-and-paste process also to twin $J14$-$10$. Now it can easily be easily seen from the examples $J14$-${10}$-$C1$ and $J14$-${10}$-$C2$ in Fig.~\ref{fig: J14_twin_cuttings} that cutting along a skew edge of two strips being twins results in the same strips of type $J15$ and $J16$. Clearly, we can draw the following 
\begin{eqnarray}
	& & \hspace*{-0.7cm}\mbox{\rm\bf{Conclusion:}} 	
	\label{lab: Conclusion_Cut_Paste} \\		%concl 
	& & \hspace*{-0.7cm}\mbox{\rm\bf{The cut-and-paste process applied to each of the twin pair strips }} \nonumber\\
	& & \hspace*{-0.7cm}\mbox{\rm\bf{of type $J14$ in Fig.~\ref{fig: J14_twin_solutions} 	results in a pair of strips of type $J15$ and $J16$}.} \nonumber
\end{eqnarray}
\newline\newline
Consequently, we can find precisely 60 different partitions for all $J15$ as well as $J16$-strips, since we have 12 twin pairs of $J14$-strips and 5 different skew cutting edges per $J14$-strip. \newline
These 60 layouts for both $J15$ and $J16$ are shown in the next figures. These layouts are arranged in the following way.\newline
We have 12 groups of layouts, each group is headed by a $J14$-layout shown in the lhs column of Fig.~\ref{fig: J14_twin_solutions}. These 12 ``header''-layouts are labeled by an alphabetical character, ranging from $A$ to $L$.
In each group 5 twin pairs consisting of a $J15$-layout and its dual $J16$-layout are given.
\clearpage
\begin{figure}[thb]
\begin{center}
\hspace{-2cm}
\includegraphics[scale = 0.7]{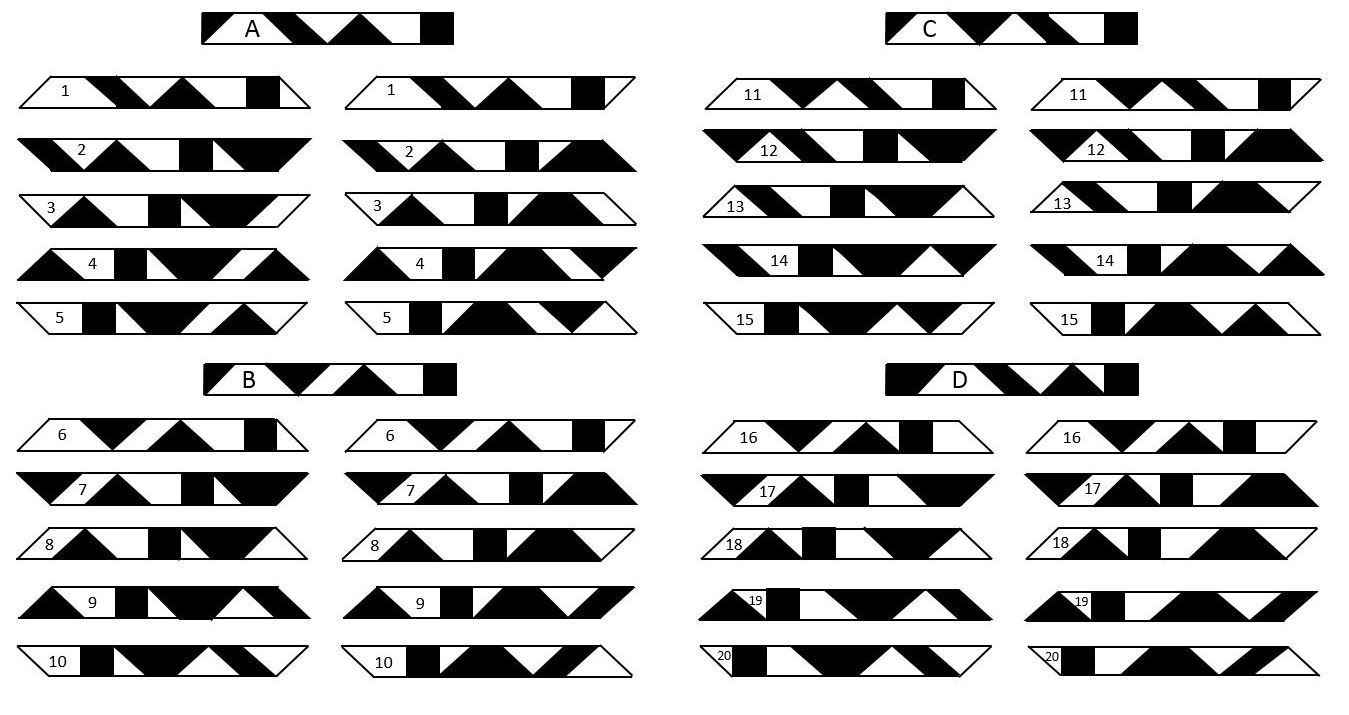}
\end{center}
\end{figure}
\begin{figure}[thb!]
\begin{center}
\hspace{-2cm}
\includegraphics[scale = 0.7]{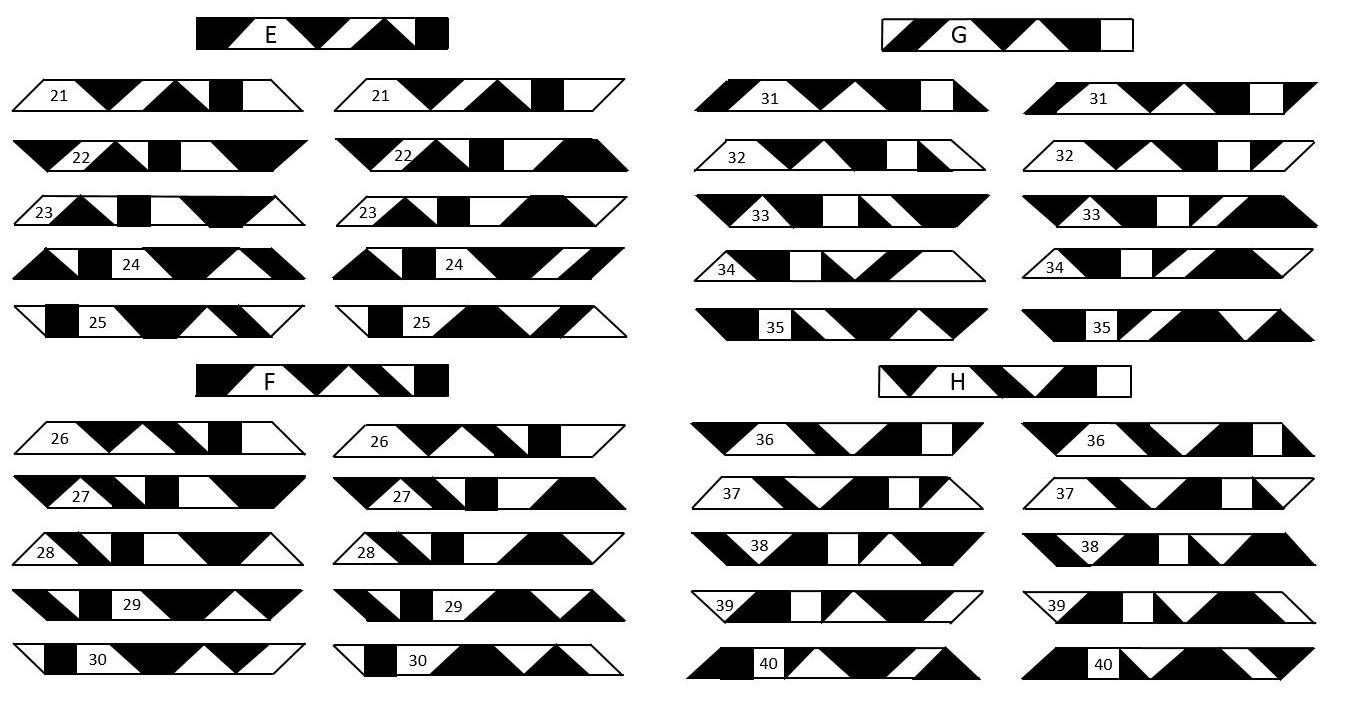}
\caption{The solutions 1-40 (out of 60) of the strips $J15$ and $J16$.}
\label{fig: J15_J16_sols_1_40}			%fig 45
\end{center}
\end{figure}
\clearpage
\begin{figure}[thb]
\begin{center}
\hspace{-2cm}
\includegraphics[scale = 0.7]{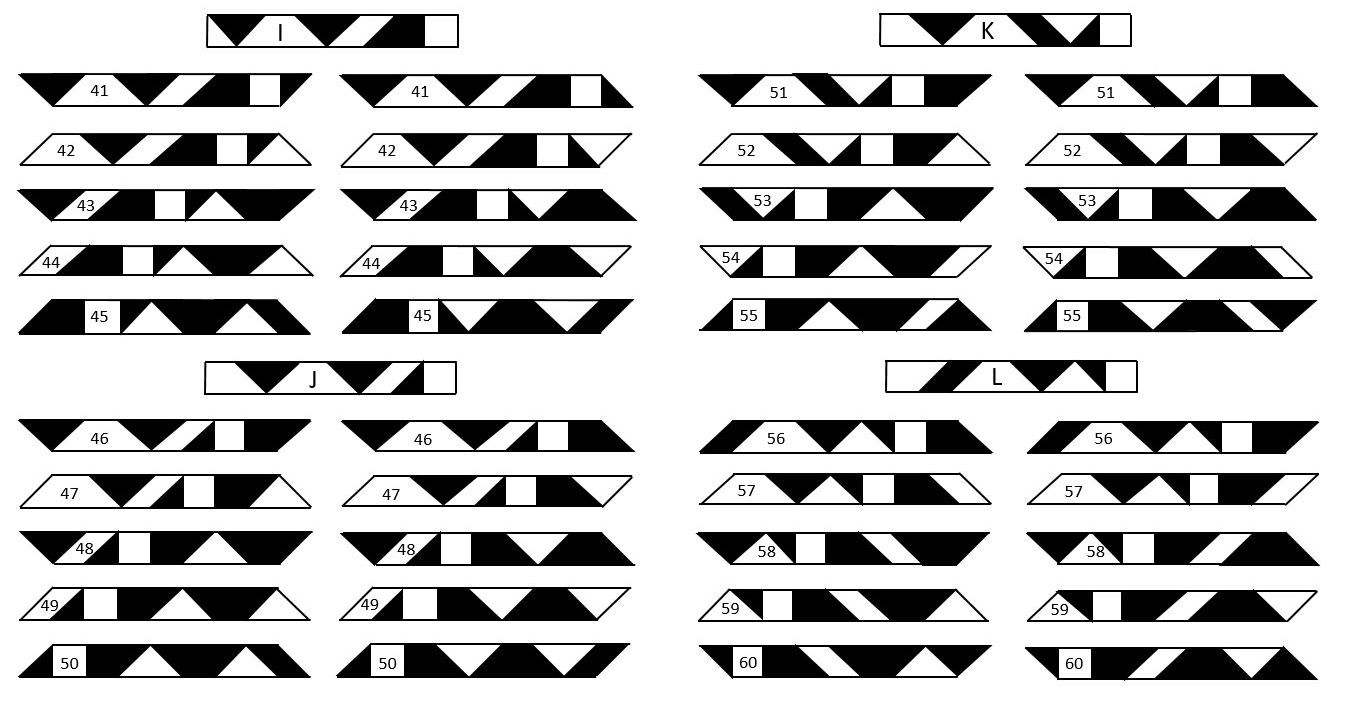}
\caption{The solutions 41-60 (out of 60) of the strips $J15$ and $J16$.}
\label{fig: J15_J16_sols_41_60}			%fig 46
\end{center}
\end{figure} 
\clearpage
\section{An algorithm for generating all partitions of a convex shape}
\label{sec: backtr_algo}
Here we will (globally) explain how the problem of finding all partitions of a convex shape can be solved by a technique used for solving a $\mathit{packing}$ problem, see \cite{Tangram:ref_Verhoeff}. This will be discussed below.
We will first start with a simple packing problem.
\subsection{A simple packing problem}
\begin{figure}[thb]		%fig in Asp_Emb.pptx
\begin{center}
\includegraphics[scale=0.45, keepaspectratio]{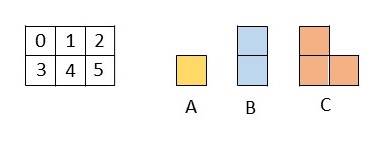}
\caption{The box with 2x3 cells and the 3 puzzle pieces $A,\; B$ and $C$.}
\label{fig: Simple_Puzzle}	%fig. 47		
\end{center}
\end{figure}

\begin{figure}[thb]		%fig in Asp_Emb.pptx
\begin{center}
\includegraphics[scale=0.4, keepaspectratio]{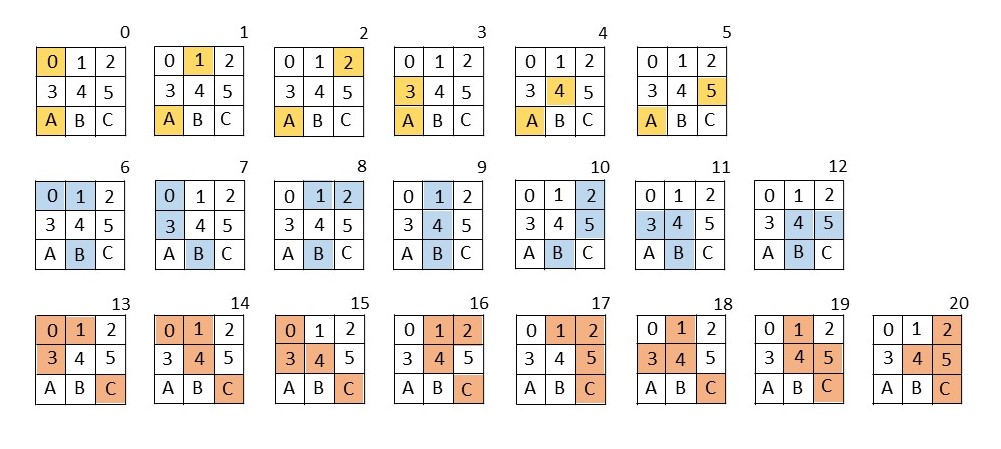}
\caption{The Aspects and Embeddings of a simple puzzle.}
\label{fig: Asp_Emb}		%fig. 48	
\end{center}
\end{figure}

Let be given a set of simple puzzle pieces and a rectangular box. The problem is to put all pieces in the box, without overlap. 
For an example, see the $Simple\; Puzzle$ in Fig.~\ref{fig: Simple_Puzzle} with a $box$ with $2x3$ unit $cells$ and the 3 $pieces$ $A,\,B$ and $C$. 
\newline
We introduce the notions $\mathit{Aspect}$ and $\mathit{Embedding}$:\newline
$\mathit{Aspect}$: \quad\quad the cell in the box, the type of a piece.\newline
$\mathit{Embedding}$: the placement of a piece in the box can be encoded by a set of $\mathit{aspects}$.\newline
In our example we have 
$Aspects\;=\; \{0\,,1,\cdots,5\,,A\,,B\,,C\}$ and the $Embeddings$ are given by Fig.~\ref{fig: Asp_Emb}. Note that a solution to such a packing puzzle consists of a set of embeddings constituting a partition of the set of aspects; that is, the embeddings in a solution are pairwise disjoint.
\begin{figure}[t!hb]
\begin{center}
\vspace{-2.cm}
\includegraphics[scale=0.48, keepaspectratio]{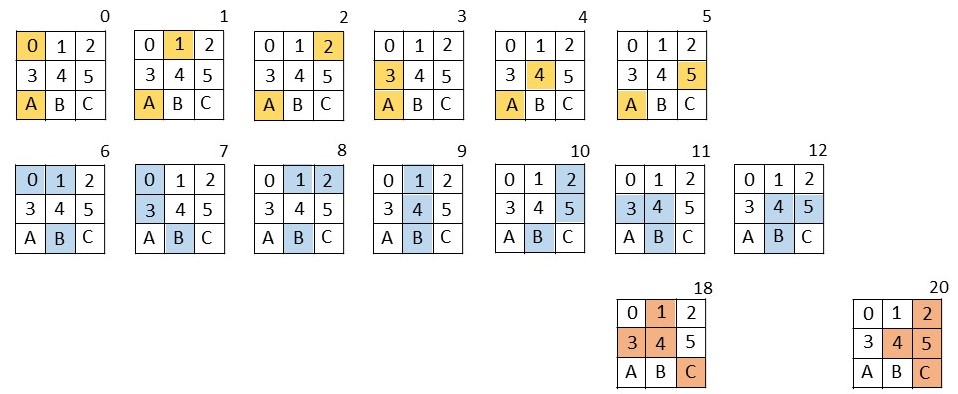}
\caption{Eliminating symmetries by restricting embeddings.}
\label{fig: Emb_restr}		%fig. 49
\end{center}
\end{figure}%

\begin{figure}[thb]
\begin{center}
\includegraphics[scale=0.48, keepaspectratio]{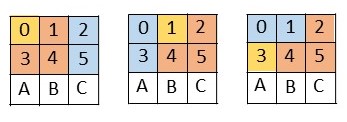}
\caption{The 3 embeddings that solve the $Simple\; L\_restricted\; Puzzle$.}
\label{fig: Sol_Simple_Puzzle}		%fig. 50
\end{center}
\end{figure}

\noindent
It is important (for performance) that we eliminate symmetries (if any). This can be done by $restricting$ the embeddings. We will illustrate this by an example. 
To this end, consider the previous example, where we will restrict piece $C$ to be the $horizontally$ $mirrored$ letter L. We call this example $Simple\; L\_restricted\; Puzzle$. This results in the embeddings in Fig.~\ref{fig: Emb_restr}.\newline
It can easily be seen that the following 3 sets of embeddings $E_1,\, E_2,\, E_3$ solve this puzzle, where 
\begin{description}
	\item $\quad\quad$ $E_1 = \{0,\; 10,\; 18 \},\;$
	$E_2 = \{1,\; 7,\; 20 \}\;$ and $\;$
	$E_3 = \{3,\; 6,\; 20 \}$. See also Fig.~\ref{fig: Sol_Simple_Puzzle}.
%	\item $\quad\quad$ See also Fig.~\ref{fig: Sol_Simple_Puzzle}.
\end{description}
\begin{figure}[thb!]
\begin{center}
\includegraphics[scale=0.5, keepaspectratio]{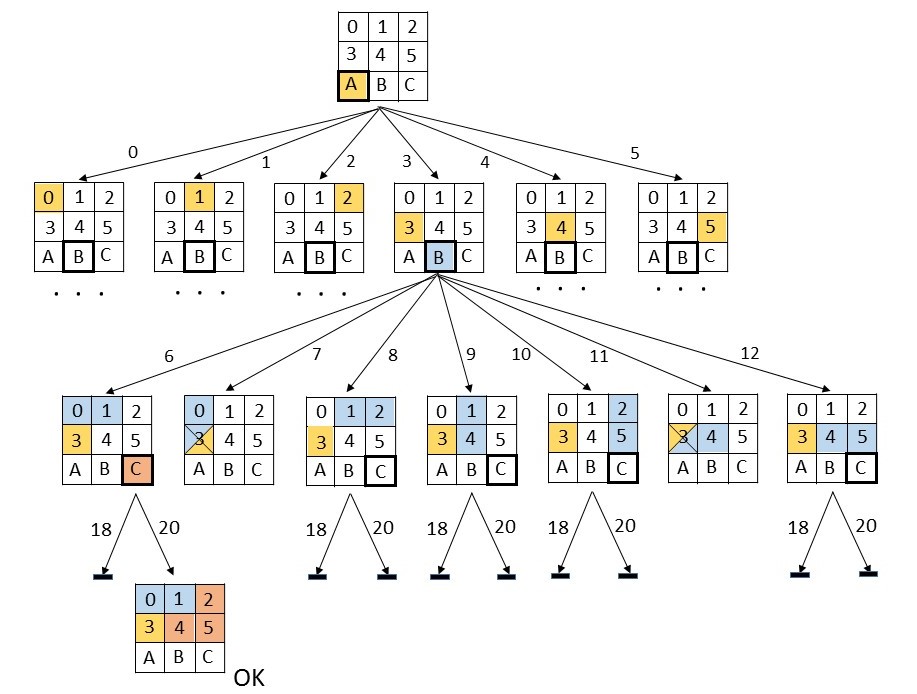}
\caption{The backtracking tree for solving the $Simple\;  L\_restricted\; Puzzle$, \newline with the embeddings \{3,\;6,\;20\} as solution.}
\label{fig: Backtracking_tree}	%fig. 51	
\end{center}
\end{figure}
\clearpage
Clearly, when dealing with a more complex puzzle we in general cannot easily find its solution(s) by hand. Then we might try to solve the puzzle using a computer and a dedicated solving procedure such as $back tracking$, see 
\cite{Tangram:ref_wikipedia_2017},
%\bibitem{Tangram:ref_wikipedia_2017}
%\url{https://en.wikipedia.org/wiki/Backtracking/}, 
%Description of the (backtracking) method, 2017,
\cite{Tangram:ref_geeks_2017}.
%\bibitem{Tangram:ref_geeks_2017}.
%\url{https://www.geeksforgeeks.org/backtracking}, 2017.
%
Recall that the $search\; tree$ is an important concept for the backtracking method. This is a graphical representation of all possible cases to be studied for finding a solution. 
%Note that the figure has the shape of an upside-down tree. 
%\newline
In Fig.~\ref{fig: Backtracking_tree} we show (a part of) the search tree corresponding to the $Simple\; L\_restricted\; Puzzle$.

\subsection{A more complicated packing problem}
\begin{figure}[thb]
\begin{center}
\includegraphics[scale=0.45, keepaspectratio]{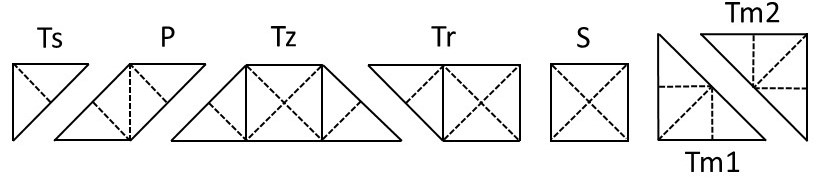}
\caption{Subdividing all Japanese tans into triangles $ts$, with $ts = Ts/2$.} 
\label{fig: Tans_triangles}		
\end{center}
\end{figure}

\begin{figure}[thb]
\begin{center}
\includegraphics[scale=0.6, keepaspectratio]{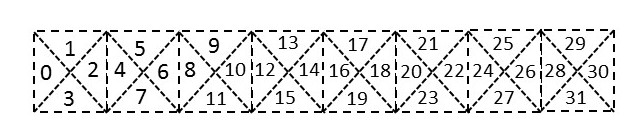}
\caption{Box with 8 cells, each being subdivided into 4 triangles $ts$.}
\label{fig: Box_with_cells}	%fig. 53	
\end{center}
\end{figure}

Now we want to consider a more complicated packing problem which is related to the problem of finding all partitions of each of the convex shapes using the Japanese set of tans.\newline
Recalling Conclusion (v) in (\ref{tan_relations}) and Fig.~\ref{fig: Chin_Jap_tans_1} in section\ref{sec: Gen_Intro} we know that each convex shape formed by the 7 Japanese tans consists of 16 isosceles rectangular triangles $Ts$. Next we can can split each $Ts$ into two smaller triangles $ts$, i.e., $ts = Ts/2$. So, the complete set of tans can be built up with 32 triangles $ts$, see Fig.~\ref{fig: Tans_triangles}. 
Now consider a rectangular box with 8 square cells, each of them being divided into 4 isosceles rectangular triangles $ts$. See Fig.~\ref{fig: Box_with_cells}.\newline
\begin{figure}[thb]
\begin{center}
\includegraphics[scale=0.55, keepaspectratio]{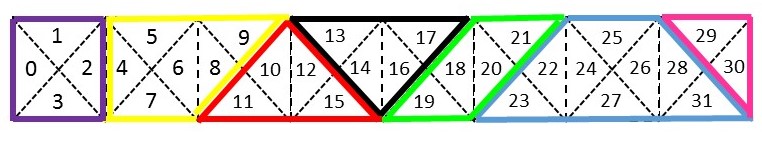}
\caption{The box in Fig.~\ref{fig: Box_with_cells} being fully covered by a partition of the strip $J14$.}
\label{fig: Tans_in_box_J14}	%fig. 54	
\end{center}
\end{figure}

\noindent
This suggests that the 7 tans can fully fill this box. Indeed, a possible filling is given in Fig.~\ref{fig: Tans_in_box_J14}.
\subsubsection{The packing problem for $Strip\; J14$}
Notice that in fact Fig.~\ref{fig: Tans_in_box_J14} shows a partition of  strip $J14$, see also Fig.~\ref{fig: All_convex_polygons}. \newline 
The situation in Fig.~\ref{fig: Tans_in_box_J14} is similar to that for the $Simple\; L\_restricted\; Puzzle$. So, similar to Fig.~\ref{fig: Backtracking_tree} for this puzzle we can also use the backtracing algorithm to find all different partitions of $J14$.
\subsubsection{The packing problem for $Strip\; J16$}
\begin{figure}[thb]
\begin{center}
\includegraphics[scale=0.6, keepaspectratio]{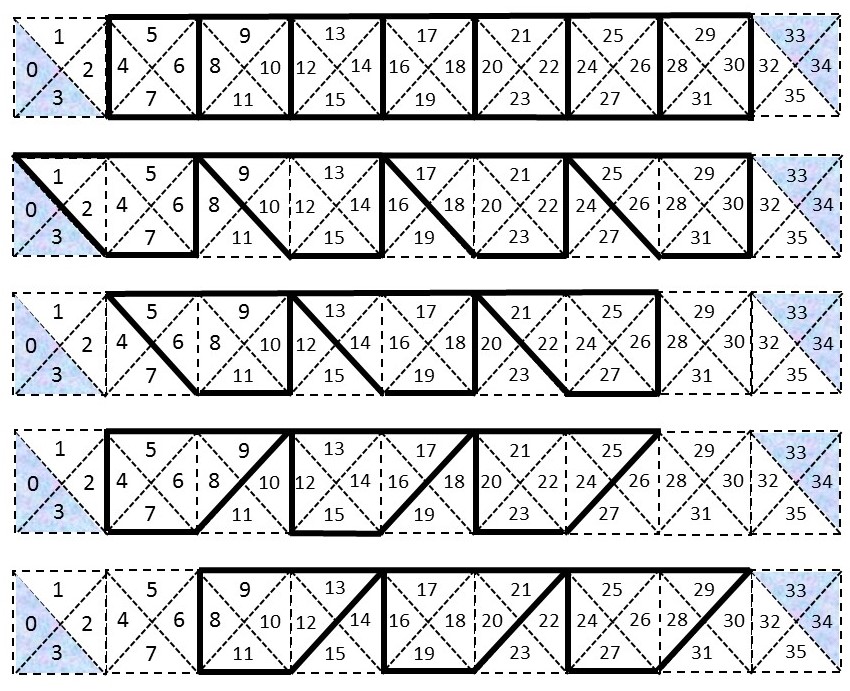}
\caption{The enlarged box being (partially) covered by a partition of the strip $J16$.}
\label{fig: Tans_in_box_J16}		
\end{center}
\end{figure}

Next, let us consider the parallelogram-shaped strip $J16$ (recall Fig.~\ref{fig: All_convex_polygons}). \newline
It is easily seen  that the partition of $J14$ in Fig.~\ref{fig: Tans_in_box_J14} can be transformed to a partition of $J16$ by moving the triangle $\{29,\,30\}$ at the most right side to the most left side. 
Of course, this partition of $J16$ fits (partially) in a $larger$ 1x9 rectangular box. See Fig.~\ref{fig: Tans_in_box_J16}.
%
% next figure in: Tans_in_box_J14_J15. pptx
%\begin{figure}[thb]
%\begin{center}
%\includegraphics[scale=0.55, keepaspectratio]{Embed_in_box_J16}
%\caption{The embeddings of $S$ and $Tr$ in box for $J16$.}
%\label{fig: Embed_in_box_J16}		
%\end{center}
%\end{figure}
\newline
Just as for the example $Simple\; L\_restricted\; Puzzle$ we can now establish the embeddings for each individual tan in the box in Fig.~\ref{fig: Tans_in_box_J16}. Notice that in order to find all essentially different partitions we have to take into account all possible orientations of each tan when establishing the embeddings. We will illustrate this for two tans ($S$ and $Tr$) only. The remaining tans can be described in a similar way.\newline
Using Fig.~\ref{fig: Tans_in_box_J16} we find the following embeddings:
\begin{itemize}
\item $\;\,S$ in row 1: $\{4,5,6,7\}, \{8,9,10,11\}, \cdots, \{28,29,30,31\}$
\item $Tr$ in row 2: $\{1,2,4,5,6,7\}, \{9,10,12,13,14,15\}, \{17,18,20,21,22,23\}, \{25,26,28,29,30,31\}$
\item $Tr$ in row 3: $\{5,6,8,9,10,11\}, \{13,14,16,17,18,19\}, \{21,22,24,25,26,27\}$ 
\item $Tr$ in row 4: $\{4,5,6,7,8,9\}, \{12,13,14,15,16,17\}, \{20,21,22,23,24,25\}$
\item $Tr$ in row 5: $\{8,9,10,11,12,13\}, \{16,17,18,19,20,21\}, \{24,25,26,27,28,29\}.$
\end{itemize}
Once having established all embeddings for all tans, we determine all possible partitions of the strip $J16$ by using a backtracking algorithm, completely similar to example $Simple\;  L\_restricted\; Puzzle$. \newline\newline
$\mathbf{Remarks:}$
\begin{itemize}
\item The establishment of the embeddings and the backtracking can be  automatically generated by a dedicated computer program.
\item Clearly, we can use the approach described above for  $all$ 16 convex polygons in Fig.~\ref{fig: All_Chin_Jap_polygons}-Right that can be formed by the set of the Japanese tans. Of course, we need a box of $mxn$ square cells for the polygons $J1$ up to $J13$, with $m$ and $n$ such that the polygon under study fits in this box.
\item In case we have two tans with the same shape but with different colour, then interchanging them in a layout gives a different partition, in contrast to the monochromatic case. This situation needs special attention when using a computer program for finding all different partitions.
\end{itemize}
In the next section we show all possible partitions of the mentioned 16 convex polygons.

\clearpage
%						pictures of all partitions
\section{All different partitions of the shapes $J01$ up to $J16$}
\begin{figure}[thb]
\begin{center}
%\hspace{-2.4cm}
%\includegraphics[scale = 0.7, angle=90]
\includegraphics[scale = 0.5]{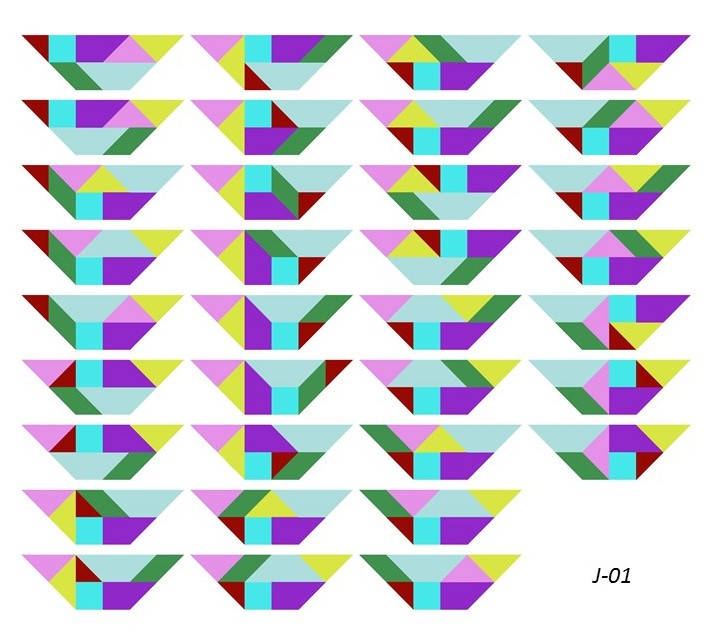}
\caption{All 34 different layouts of shape $J01$.}
\label{fig: shape-J01}		% shape J-01
\end{center}
\end{figure}
\begin{figure}[thb!]
\begin{center}
%\hspace{-2.4cm}
%\includegraphics[scale = 0.7, angle=90]
\includegraphics[scale = 0.45]{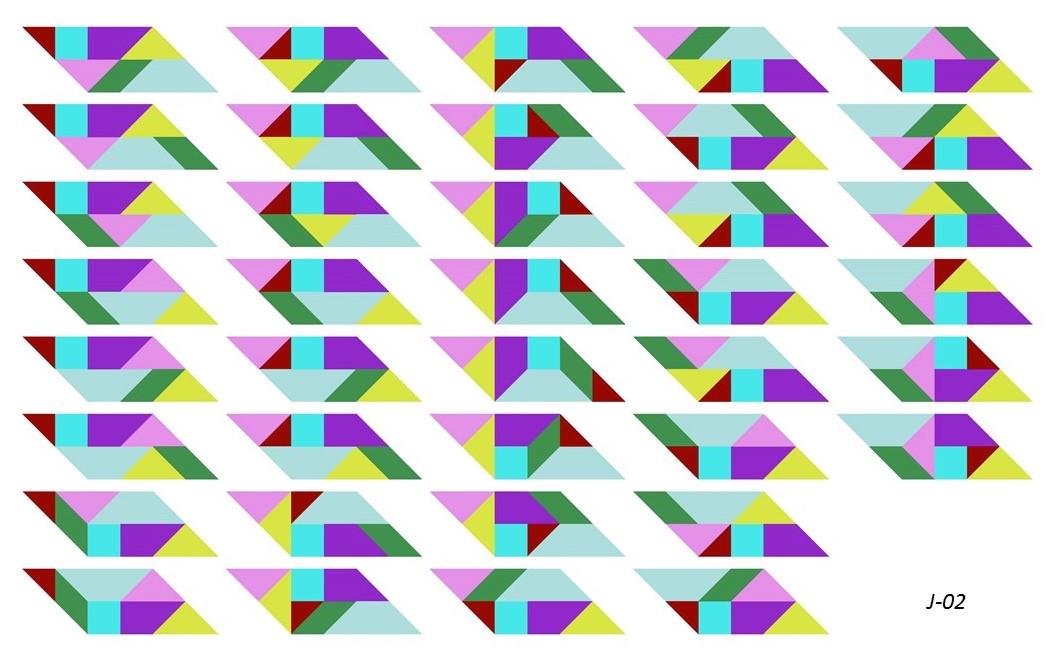}
\caption{All 38 different layouts of shape $J02$.}
\label{fig: shape-J02}		% shape J-02
\end{center}
\end{figure}
\clearpage
\begin{figure}[thb]
\begin{center}
%\hspace{-2.4cm}
%\includegraphics[scale = 0.7, angle=90]
\includegraphics[scale = 0.4]{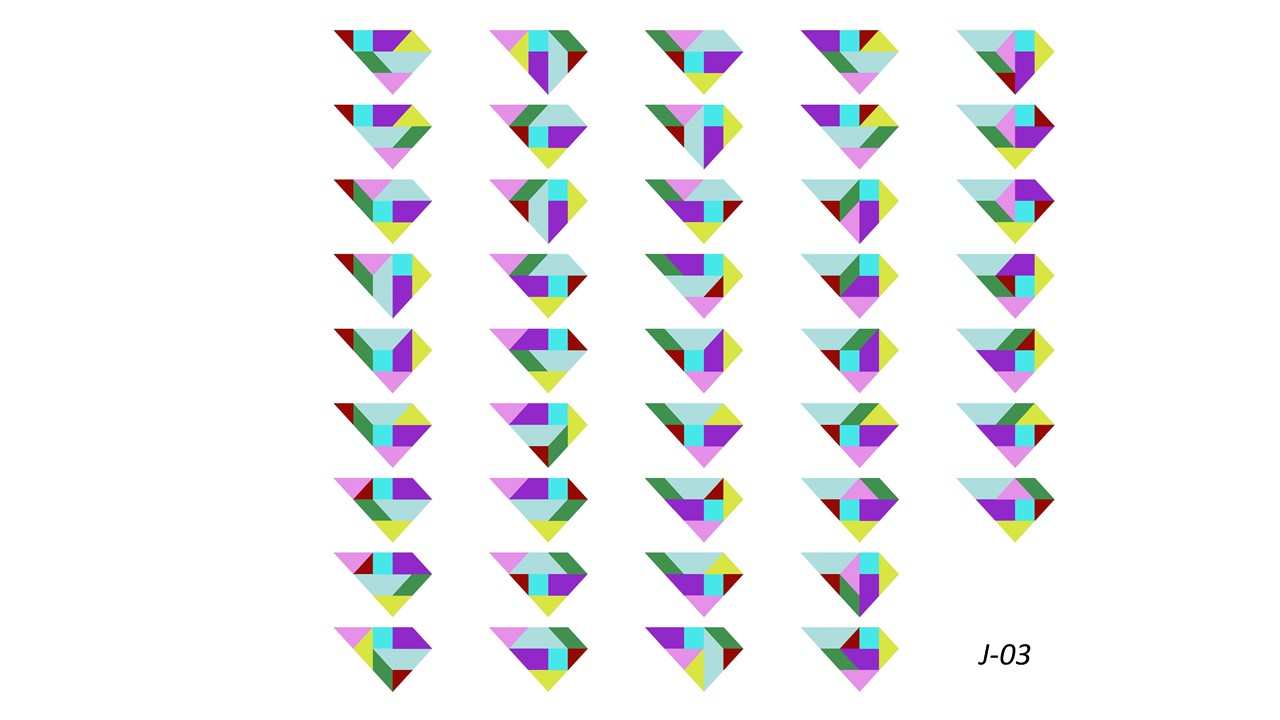}
\caption{All 43 different layouts of shape $J03$.}
\label{fig: shape-J03}		% shape J-03
\end{center}
\end{figure}
\begin{figure}[thb!]
\begin{center}
%\hspace{-2.4cm}
%\includegraphics[scale = 0.7, angle=90]
\includegraphics[scale = 0.52]{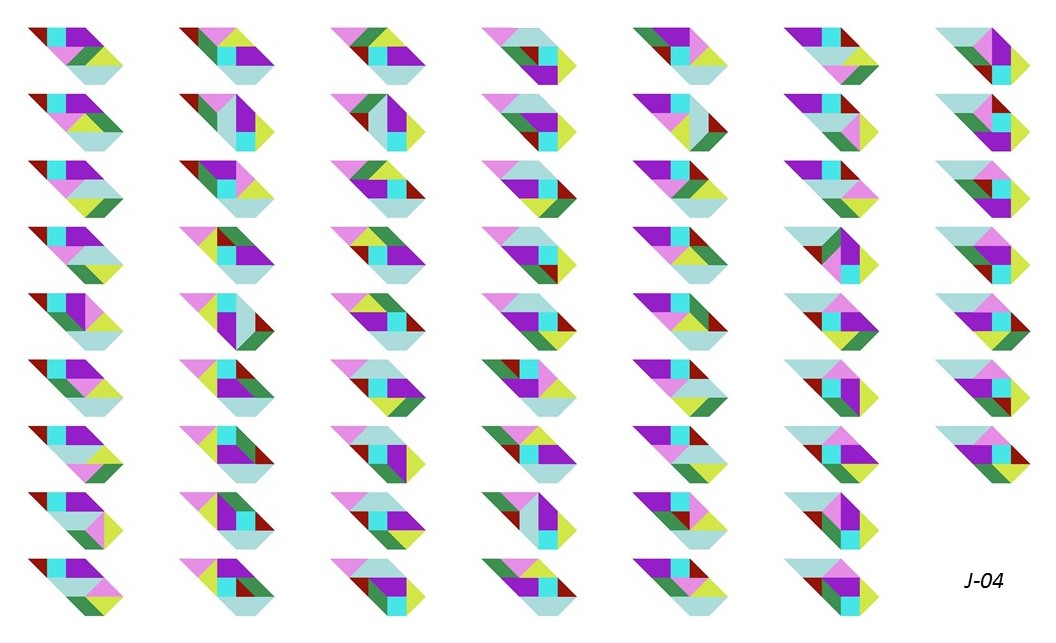}
\caption{All 61 different layouts of shape $J04$.}
\label{fig: shape-J04}		% shape J-04
\end{center}
\end{figure}
\clearpage
\begin{figure}[thb]
\begin{center}
%\hspace{-2.4cm}
%\includegraphics[scale = 0.7, angle=90]
\includegraphics[scale = 0.4]{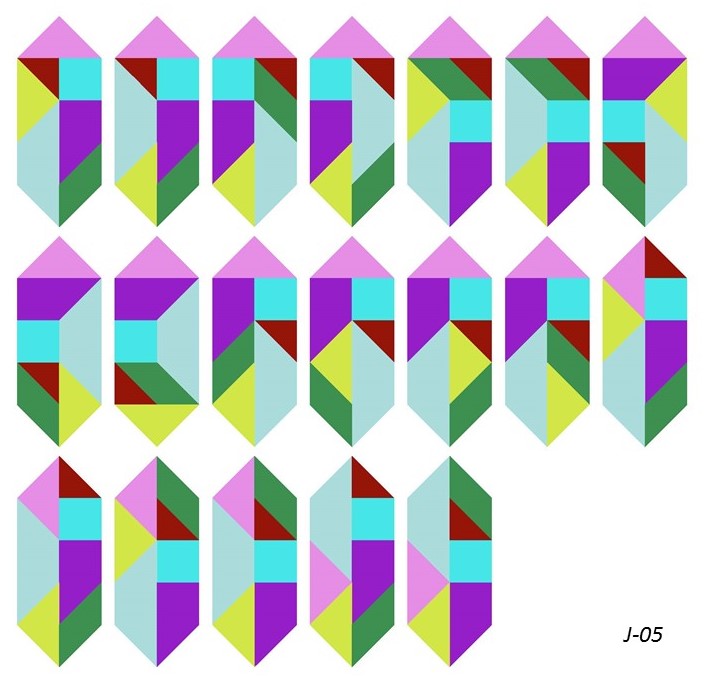}
\caption{All 19 different layouts of shape $J05$.}
\label{fig: shape-J05}		% shape J-05
\end{center}
\end{figure}
\begin{figure}[thb!]
\begin{center}
\hspace{-1.9cm}
\includegraphics[scale = 0.58]{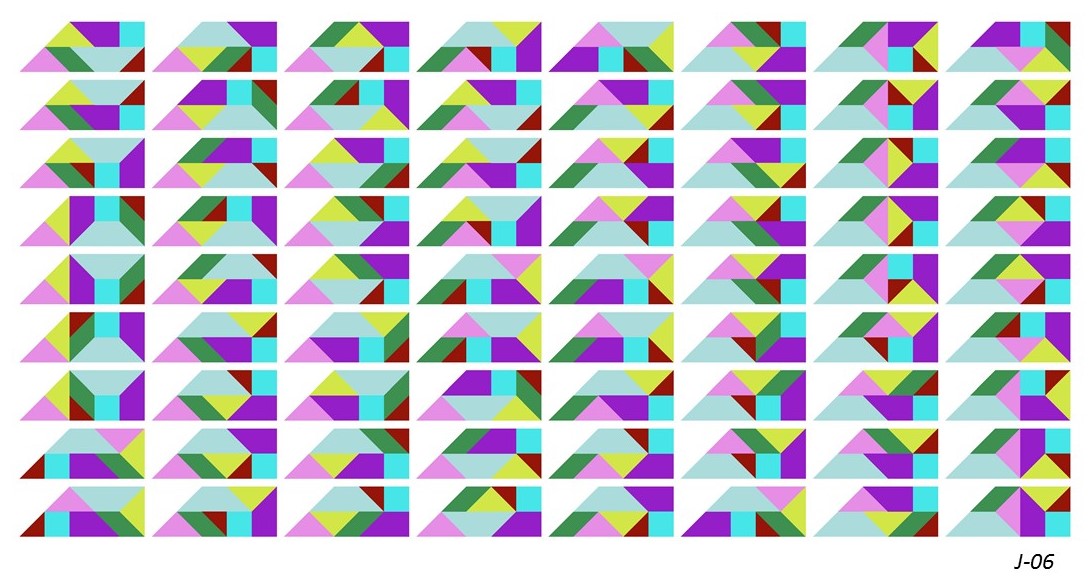}
\caption{All 72 different layouts of shape $J06$.}
\label{fig: shape-J06}		% shape J-06
\end{center}
\end{figure}
\clearpage
\begin{figure}[thb]
\begin{center}
%\hspace{-2.4cm}
%\includegraphics[scale = 0.7, angle=90]
\includegraphics[scale = 0.5]{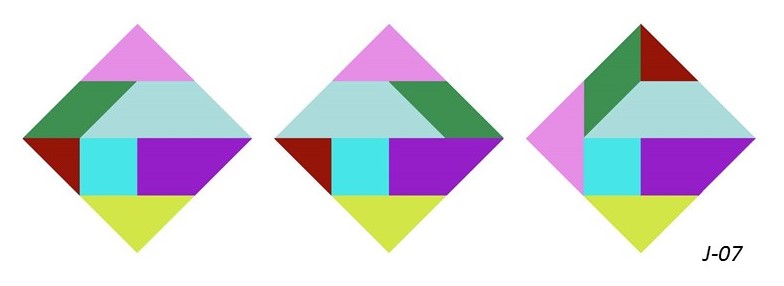}
\caption{All 3 different layouts of shape $J07$.}
\label{fig: shape-J07}		% shape J-07
\end{center}
\end{figure}
\begin{figure}[thb!]
\begin{center}
\hspace{-1.5cm}
\includegraphics[scale = 0.58]{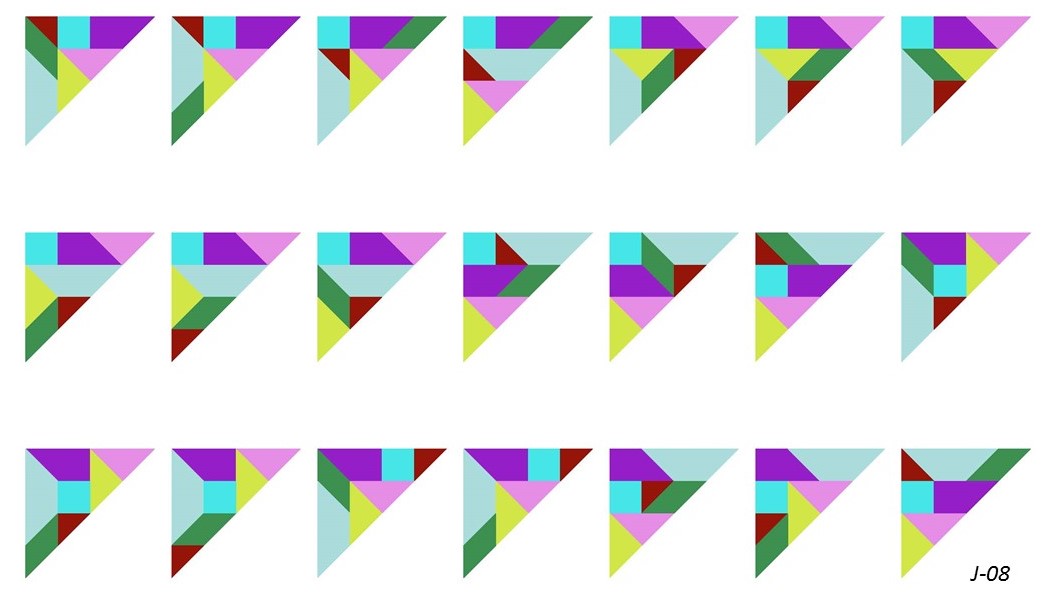}
\caption{All 21 different layouts of shape $J08$.}
\label{fig: shape-J08}		% shape J-08
\end{center}
\end{figure}
\clearpage
\begin{figure}[thb]
\begin{center}
%\hspace{-2.4cm}
%\includegraphics[scale = 0.7, angle=90]
\includegraphics[scale = 0.47]{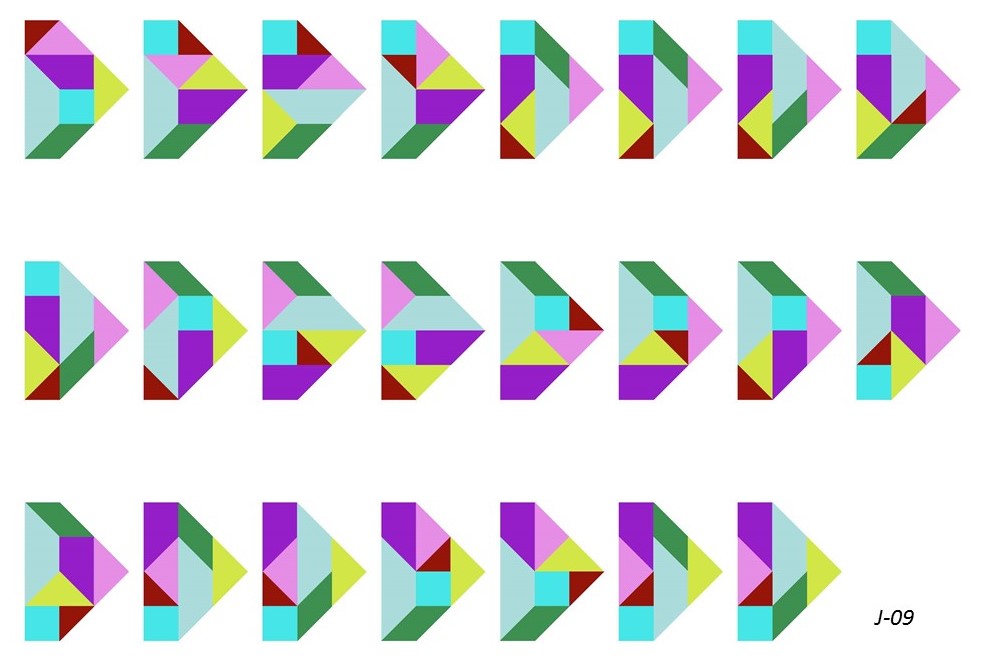}
\caption{All 23 different layouts of shape $J09$.}
\label{fig: shape-J09}		% shape J-09
\end{center}
\end{figure}
\begin{figure}[thb!]
\begin{center}
%\hspace{-1.5cm}
%\includegraphics[scale = 0.7, angle=90]
\includegraphics[scale = 0.47]{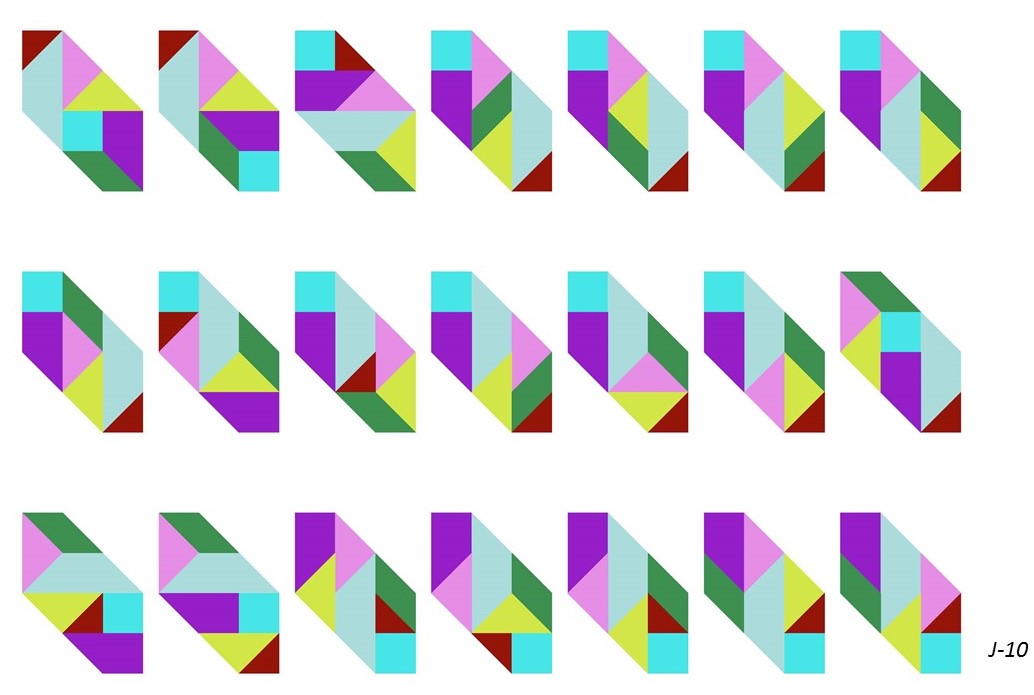}
\caption{All 21 different layouts of shape $J10$.}
\label{fig: shape-J10}		% shape J-10
\end{center}
\end{figure}
\clearpage
\begin{figure}[thb]
\begin{center}
%\hspace{-2.4cm}
%\includegraphics[scale = 0.7, angle=90]
\includegraphics[scale = 0.5]{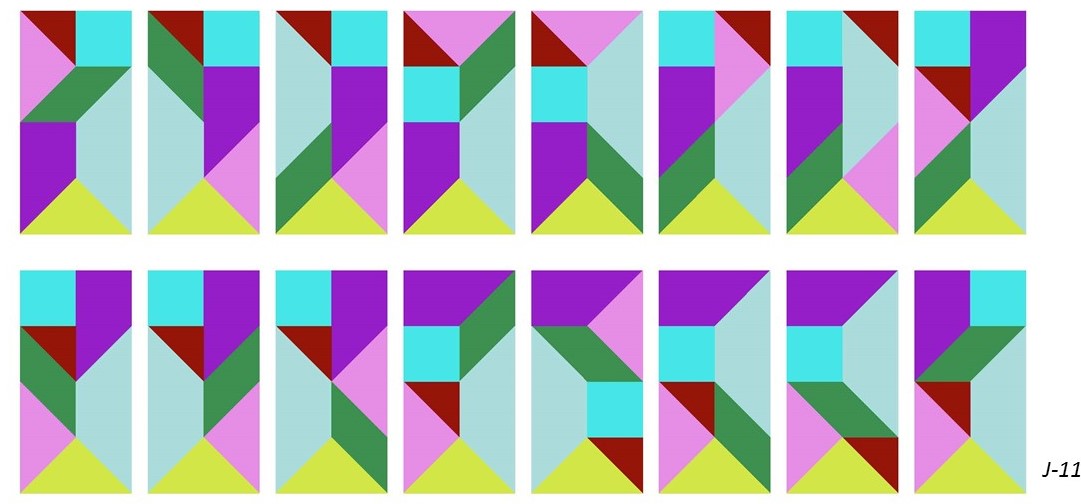}
\caption{All 16 different layouts of shape $J11$.}
\label{fig: shape-J11}		% shape J-11
\end{center}
\end{figure}
\begin{figure}[thb]
\begin{center}
%\hspace{-1.5cm}
%\includegraphics[scale = 0.7, angle=90]
\includegraphics[scale = 0.5]{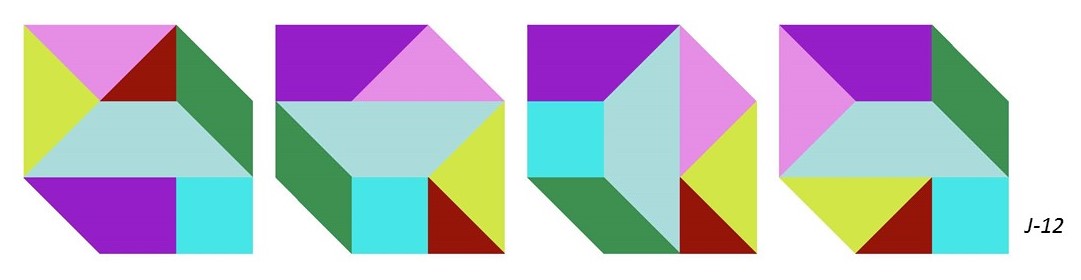}
\caption{All 4 different layouts of shape $J12$.}
\label{fig: shape-J12}		% shape J-12
\end{center}
\end{figure}
\clearpage
\begin{figure}[thb!]
\begin{center}
\hspace{-1.8cm}
\includegraphics[scale = 0.62]{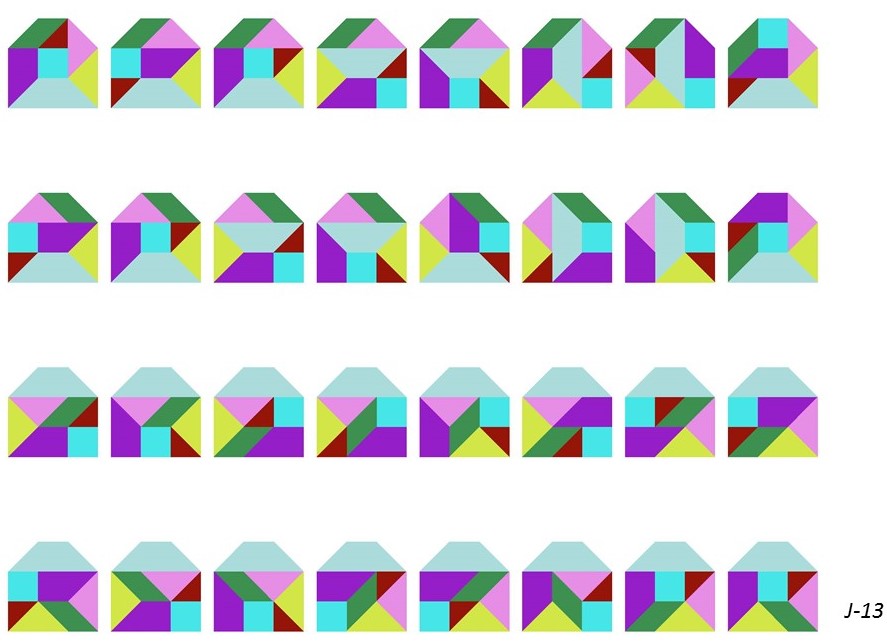}
\caption{All 32 different layouts of shape $J13$.}
\label{fig: shape-J13}		% shape J-13
\end{center}
\end{figure}
\begin{figure}[thb!]
\begin{center}
\hspace{-2cm}
\includegraphics[scale = 0.6]{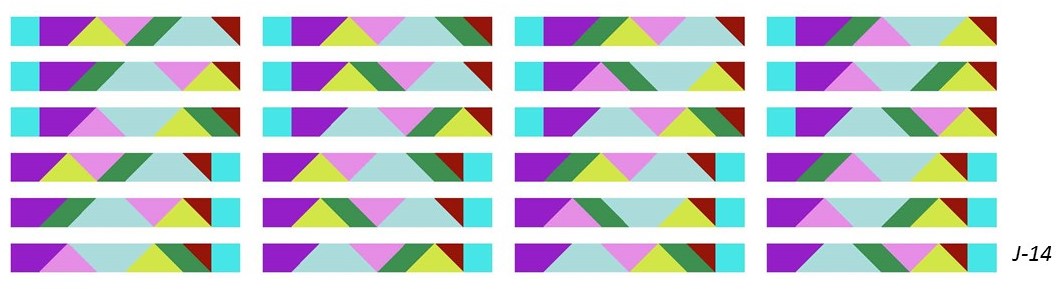}
\caption{All 24 different layouts of shape $J14$.}
\label{fig: shape-J14}		% shape J-14
\end{center}
\end{figure}
\clearpage
\begin{figure}[thb]
\begin{center}
%\hspace{-1.5cm}
%\includegraphics[scale = 0.7, angle=90]
\includegraphics[scale = 0.478]{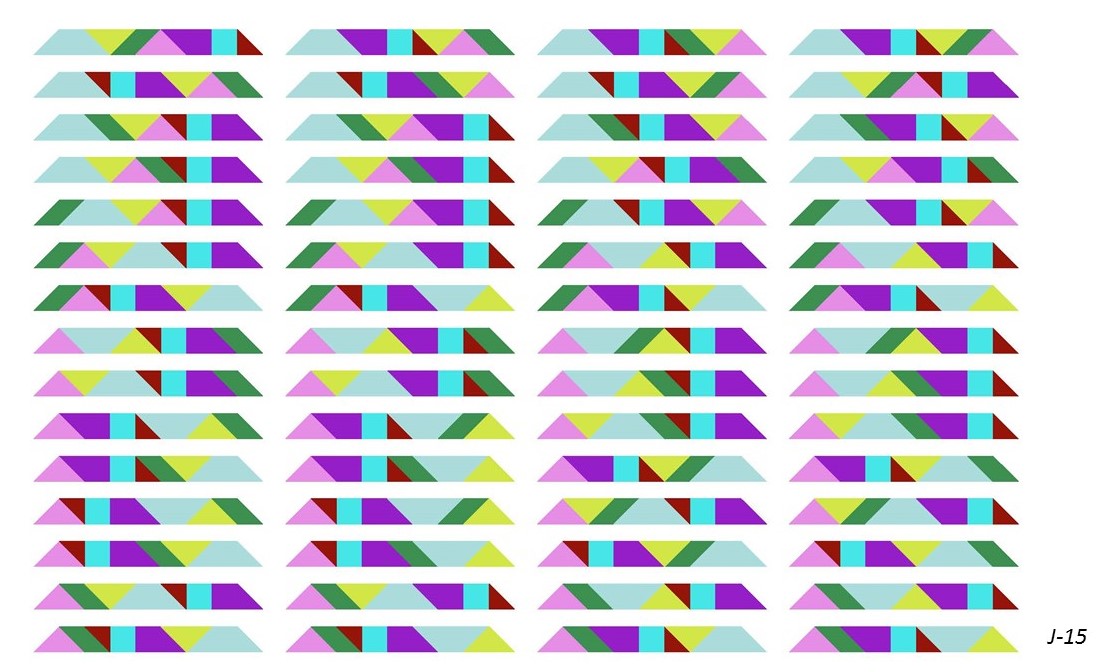}
\caption{All 60 different layouts of shape $J15$.}
\label{fig: shape-J15}		% shape J-15
\end{center}
\end{figure}
\begin{figure}[thb!]
\begin{center}
%\hspace{-1.5cm}
%\includegraphics[scale = 0.7, angle=90]
\includegraphics[scale = 0.48]{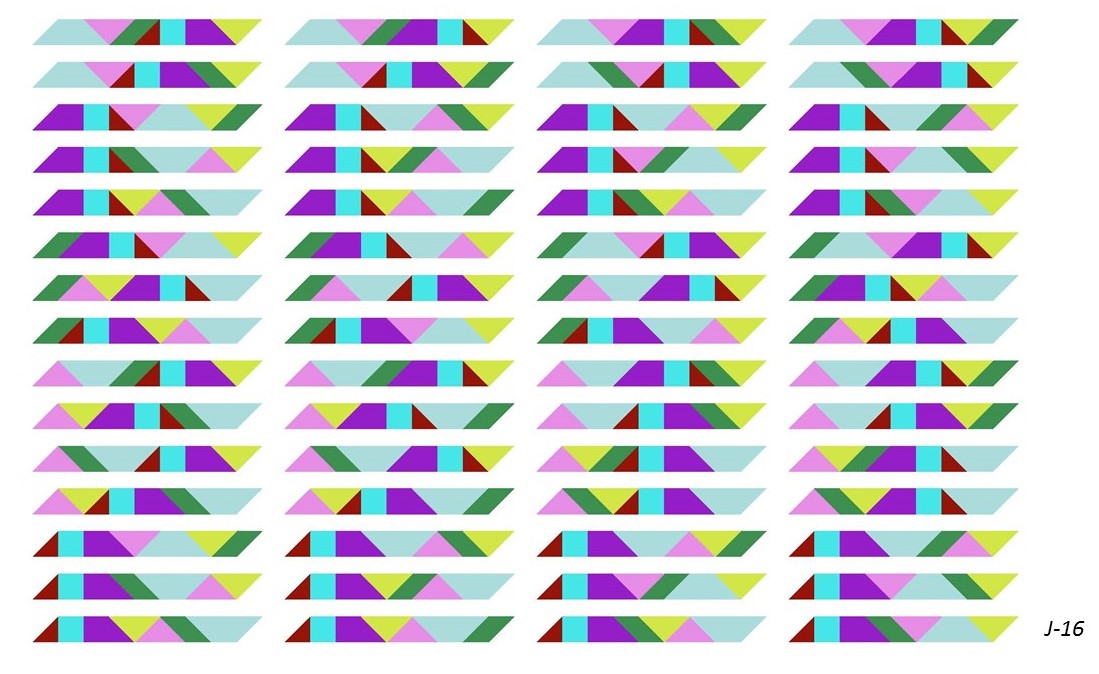}
\caption{All 60 different layouts of shape $J16$.}
\label{fig: shape-J16}		% shape J-16
\end{center}
\end{figure}
\clearpage
\begin{table}[thb!]
\hspace{-2cm}
\caption{Summary: Number of different partitions per shape.}
\vspace{2mm}
\centering
\begin{tabular}{|c |c c c c c c c c c c c c c c c c|}
\hline
	&  &   &   &    &    &    & 	&   &    &    &    &  &  &  &  &	\\ %empty line
%	     &  &   &   &    &    &    & 	&   &    &    &    &  &  &  &  &		\\
	$Shape$ 
	&$1$  &$2$  &$3$  &$4$   &$5$   &$6$ &$7$ &$8$
						   &$9$  &$10$ &$11$ &$12$  
						   &$13$  &$14$ &$15$ &$16$\\
\hline
	    &  &   &   &    &    &    & 	&   &    &    &    &  &  &  &  &  \\
%empty line
% $Number\; of$  &  &   &   &    &    &    & 	&   &    &    &    &  &  &  &  & \\
 \# $partitions$ 
 &34 &38 &43 &61 &19 &72 &3 &21 &23 &21 &16 &4 &32 &24 &60 &60\\
\hline
\end{tabular}
\label{Count_diff_sols_per_shape}
\end{table}
\clearpage
\section{Acknowledgments and References}\label{sec: Refs}
\mbox{\rm\bf\Large{Acknowledgments}}
\newline\newline
The author wishes to thank co-author Tom Verhoeff for his stimulating discussions and for providing and running his dedicated software to generate all solutions we were looking for. Without his help the results as presented in this report would not have been found.\newline
Last but not least I would like to express my appreciation to Mrs. Hui Jun Chang who read this report carefully and gave several valuable suggestions for improvement. 
%
%
%
 % use this if you have 1-9 references


\begin{thebibliography}{99} % use this if you have 1-9 references
% \begin{thebibliography}{10} % use this if you have 10-99 references
%\item[1]
\bibitem{Tangram:ref_paper_1942}
Fu~Traing~Wang, Chuan-Chih~Hsiung, \emph{A Theorem on the Tangram}, The American Mathematical Monthly, Vol. 49, No. 9, 1942, pp. 596-599, available at http://www.jstor.org/stable/2303340.
%
%\item[2]
\bibitem{Tangram:ref_elffers_1976}
Joost~Elffers, \emph{Tangram, Das alte chinesische Formenspiel / Het oude Chinese vormenspel}, M.~DuMont Schauberg, Cologne, 1976, 3rd edition.
%
%\item[3]
\bibitem{Tangram:ref_beelen_2017}
T.G.J. Beelen, \emph{Finding all convex tangrams}, Eindhoven University of Technology, CASA-report 17-02, 2017.
%
%\item[4]
\bibitem{Tangram:ref_beelen_2018}
T.G.J. Beelen, \emph{Determining the essentially different partitons of all Chinese convex tangrams}, Eindhoven University of Technology, CASA-report 18-02, 2018.
%\item[4b]
\bibitem{Tangram:ref_beelen_verhoeff_2018}
T.G.J. Beelen, T. Verhoeff, \emph{Determining the essentially different partitons of all Japanese convex tangrams}, Eindhoven University of Technology, CASA-report 18-07, 2018.
%
%\item[5]
\bibitem{Tangram:ref_pentoma_2017}
\url{http://www.pentoma.de/}, Tangram Figuren Uebersicht, Konvexe Figuren, 2017.
%
%\item[6]
\bibitem{Tangram:ref_Verhoeff}
Erik van der Tol and Tom Verhoeff, \emph{The Puzzle Processor Project: Towards an Implementation}, Nat. Lab. Unclassified Report NL-UR 2000/828, Koninklijke Philips Electronics N.V., 2001.
%Erik van der Tol and Tom Verhoeff. The Puzzle Processor Project: Towards an Implementation. Nat. Lab. Unclassified Report NL-UR 2000/828, Koninklijke Philips Electronics N.B., 2001.
%
%\item[7]
\bibitem{Tangram:ref_wikipedia_2017}
%https:$/$en.wikipedia.org$/$wiki$/$Backtracking$/$Description_of_the_method
\url{https://en.wikipedia.org/wiki/Backtracking/}, Description of the (backtracking) method, 2017.
%\url{https://en.wikipedia.org/wiki/Backtracking/}, Description of the (backtracking) method, 2017.
%
%\item[8]
\bibitem{Tangram:ref_geeks_2017}
\url{https://www.geeksforgeeks.org/backtracking}, 2017.
%
%\item[9]
\bibitem{Tangram:ref_epstein_uehara_2014}
Eli~Fox-Epstein, Ryuhei~Uehara, \emph{The Convex Configurations of ``Sei Shonagon Chie no Ita''}, and Other Dissection Puzzles, arXiv:1407.1923v1, July 2014.
%%
%\item[10]
%Jerry~Slocum, \emph{The Tangram Book: The Story of the Chinese Puzzle over 2000 Puzzles to Solve}, Sterling Publishing, 2004.
%
%\end{description}
%
%
\end{thebibliography}
\end{document}